\documentclass[a4paper,11pt,reqno]{amsart}

\newtheorem{theorem}[equation]{Theorem}

\newtheorem{corollary}[equation]{Corollary}
\newtheorem{proposition}[equation]{Proposition}
\newtheorem{conjecture}[equation]{Conjecture}
\numberwithin{equation}{section}

\theoremstyle{definition}
\newtheorem{definition}[equation]{Definition}

\newtheorem{examples}[equation]{Examples}
\newtheorem{remark}[equation]{Remark}

\newcommand{\bZ}{{\mathbb Z}}

\newcommand{\frg}{{\mathfrak g}}
\newcommand{\frtg}{{\tilde{\mathfrak g}}}

\newcommand{\frs}{{\mathfrak s}}
\newcommand{\fra}{{\mathfrak a}}
\newcommand{\frz}{{\mathfrak z}}
\newcommand{\frL}{{\mathfrak L}}
\newcommand{\frgo} {{\frg_{\bar 0}}}
\newcommand{\frguno} {{\frg_{\bar 1}}}

\newcommand{\calJord}{{\mathcal{J}ord}}
\newcommand{\subo}{_{\bar 0}}
\newcommand{\subuno}{_{\bar 1}}

\providecommand{\espan}[1]{\text{span}\left\{ #1\right\}}

 \DeclareMathOperator{\frsl}{{\mathfrak{sl}}}
 \DeclareMathOperator{\frsp}{{\mathfrak{sp}}}
 \DeclareMathOperator{\frso}{{\mathfrak{so}}}
 \DeclareMathOperator{\frpsl}{{\mathfrak{psl}}}
 \DeclareMathOperator{\frgl}{{\mathfrak{gl}}}
 \DeclareMathOperator{\frpgl}{{\mathfrak{pgl}}}
 \DeclareMathOperator{\frosp}{{\mathfrak{osp}}}
 
 \DeclareMathOperator{\fru}{{\mathfrak{u}}}

 \DeclareMathOperator{\trace}{trace}
 \DeclareMathOperator{\ad}{ad}
 \DeclareMathOperator{\der}{der}
 \DeclareMathOperator{\inder}{inder}
 \DeclareMathOperator{\End}{End}
 \DeclareMathOperator{\Hom}{Hom}
 \DeclareMathOperator{\Mat}{Mat}

 \DeclareMathOperator{\charac}{char}

\newcommand{\T}{\bigl(T,[...],(.\vert.)\bigr)}
\newcommand{\STS}{symplectic triple system}
\newcommand{\STSs}{symplectic triple systems}
\newcommand{\OTS}{orthogonal triple system}
\newcommand{\OTSs}{orthogonal triple systems}

\newenvironment{romanenumerate}
 {\begin{enumerate}
 }{\end{enumerate}}

\newenvironment{alphaenumerate}
 {\begin{enumerate}
 }{\end{enumerate}}

\begin{document}


\title{New simple Lie superalgebras in characteristic $3$}

\author{Alberto Elduque}
 \thanks{Supported by the Spanish Ministerio de Educaci\'{o}n y Ciencia
 and FEDER (MTM 2004-081159-C04-02 and ) and by the
Diputaci\'on General de Arag\'on (Grupo de Investigaci\'on de
\'Algebra)}
 \address{Departamento de Matem\'aticas, Universidad de
Zaragoza, 50009 Zaragoza, Spain}
 \email{elduque@unizar.es}

\date{\today}

\subjclass[2000]{Primary 17B50; Secondary 17B60, 17B25}

\keywords{Simple, Lie, algebra, superalgebra, symplectic triple
system, orthogonal triple system, Freudenthal}

\begin{abstract}
Symplectic (respectively orthogonal) triple systems provide
constructions of Lie algebras (resp. superalgebras). However, in
characteristic $3$, it is shown that this role can be interchanged
and that Lie superalgebras (resp. algebras) can be built out of
\STSs\ (resp. \OTSs) with a different construction. As a
consequence, new simple finite dimensional Lie superalgebras, as
well as new models of some nonclassical simple Lie algebras, over
fields of characteristic $3$, will be obtained.
\end{abstract}

\maketitle


\section{Introduction}

Symplectic triple systems \cite{YamAs} are basic ingredients in
the construction of some $5$-graded Lie algebras. They consist of
a vector space $T$ endowed with a trilinear product $[...]$ and a
nonzero alternating bilinear form $(.\vert.)$ satisfying some
conditions (Definition \ref{df:STS}). Given such a system over a
field $k$ of characteristic $\ne 2$, a Lie algebra can be defined
on
\begin{equation}\label{eq:introg}
\frg=\bigl(\frsl_2(k)\oplus \inder(T)\bigr)\oplus (k^2\otimes T),
\end{equation}
where $\inder(T)=\espan{[xy.]: x,y\in T}$ and $k^2$ is the natural
module for the three dimensional simple split Lie algebra
$\frsl_2(k)$ (Theorem \ref{th:STSLie}). The $5$-grading is
obtained by looking at the eigenspaces of the adjoint action of a
Cartan subalgebra in $\frsl_2(k)$. Equivalent formulations of this
construction have been given by means of Freudenthal triple
systems \cite{FrII,FrVIII} or some other ternary algebras
\cite{Faulkner}. In this way, different models of the exceptional
simple Lie algebras have been obtained. (For constructions based
on different systems see, for instance, \cite{Kantor} or
\cite{Allison}.)

Here, besides the classical exceptional simple Lie algebras, it
will be shown that, over fields of characteristic $3$, a
one-parameter family of ten dimensional simple Lie algebras
discovered by Kostrikin \cite{Kostrikin}, as well as a simple Lie
algebra of dimension $29$ considered by Brown \cite{Brown29}, can
be obtained, respectively, from a family of two dimensional simple
\STSs\ (Proposition \ref{pr:Lepsilon}) and from an eight
dimensional simple \STS\ (Remark \ref{re:STSLiechar3}), thus
providing very simple descriptions of these Lie algebras.

But the main point of this paper is the description of another
feature that occurs in characteristic $3$. It will be noticed that
over fields of this characteristic, one can forget about the
$\frsl_2(k)$ and  the $k^2$ in \eqref{eq:introg} and consider
instead the direct sum
\begin{equation}\label{eq:introgtilde}
\frtg=\inder(T)\oplus T
\end{equation}
with a natural multiplication. Then $\frtg$ becomes a Lie
superalgebra (instead of a Lie algebra). In other words, any \STS\
becomes an anti-Lie triple system in characteristic $3$. For
specific symplectic triple systems, new simple Lie superalgebras
will be obtained over fields of characteristic $3$.

Dually, orthogonal triple systems (defined by Okubo \cite{OkuboI})
consist of a vector space endowed with a triple product and a
bilinear symmetric form. The defining relations of an \OTS\ are
obtained by changing symmetry and skew-symmetry in the defining
relations of a \STS\ (Definition \ref{df:OTS}). It turns out that
one can consider for \OTSs\ the same construction as in
\eqref{eq:introg} and this gives now $5$-graded Lie superalgebras
(instead of algebras) The basic classical exceptional Lie
superalgebras in characteristic $0$ are obtained in this way
\cite{KamiyaOkubo,EKO}.

Again, in characteristic $3$, some new simple Lie superalgebras
appear with this construction, which are associated to a family of
\OTSs\ related to central simple Jordan algebras of degree $3$
(Theorem \ref{th:othernewchar3}).

And, as for \STSs,  in characteristic $3$ one can forget about
$\frsl_2(k)$ and $k^2$ in \eqref{eq:introg} and get $\frtg$ as in
\eqref{eq:introgtilde} which is, now, a Lie algebra (instead of a
superalgebra). New constructions of some classical and
nonclassical simple Lie algebras are obtained in this way
(Examples \ref{ex:OTSexamples}).

\smallskip

\emph{Throughout the paper, all the vector spaces and algebras
will be assumed to be finite dimensional over a ground field $k$
of characteristic $\ne 2$.}

\smallskip

In Section \ref{se:STS}, the definition of \STSs\ will be
recalled, as well as the construction of the Lie algebras $\frg$
in \eqref{eq:introg}. The relationship of \STSs\ with Freudenthal
triple systems and with a class of ternary algebras defined by
Faulkner will be given, and this will lead to the classification
of the simple \STSs\ over algebraically closed fields. In
characteristic $3$, there appears a new family of two dimensional
\STSs\ and a new eight dimensional simple \STS, which give rise,
through the construction in \eqref{eq:introg}, to the
one-parameter family of Kostrikin simple Lie algebras and to the
$29$ dimensional simple Lie algebra of Brown.

Section \ref{se:STSLie} is devoted to the construction
\eqref{eq:introgtilde} of Lie superalgebras over fields of
characteristic $3$. New simple Lie superalgebras of dimension
$18$, $35$, $54$, $98$ and $189$ appear with this construction.
Besides, one can relax the definition of a \STS\ by allowing the
alternating bilinear form to be trivial. It will be shown that no
new simple \STSs\ appear in characteristic $\ne 3$, and it is
conjectured that only two dimensional simple systems are possible
in characteristic $3$.

Section \ref{se:OTS} is devoted to \OTSs\ and the construction
\eqref{eq:introg} of Lie superalgebras from them. The
classification of the simple $(-1,-1)$ balanced Freudenthal Kantor
triple systems in \cite{EKO} immediately provides the
classification of the simple \OTSs\ over fields of characteristic
$0$. All the simple systems that appear in this classification
have their counterparts in prime characteristic, and in
characteristic $3$ a new family of simple \OTSs\ will be defined,
in terms of central simple Jordan algebras of degree $3$. However,
there is no known classification of the simple \OTSs\ over fields
of prime characteristic. The construction \eqref{eq:introg}
provides here, in particular, new simple Lie superalgebras of
dimension $24$, $37$, $50$ and $105$ over fields of characteristic
$3$.

Finally, Section \ref{se:OTSLie} will deal with the construction
\eqref{eq:introgtilde} of Lie algebras in term of \OTSs. Apart
from some classical Lie algebras, new models of Kostrikin's ten
dimensional simple Lie algebras and of Brown's $29$ dimensional
simple Lie algebra will be given.

\smallskip

In a future work, some of the new simple Lie superalgebras in
characteristic $3$ that will be constructed here, will be shown to
appear in a natural extension of Freudenthal's Magic Square,
obtained by means of composition superalgebras \cite{EO}.


\section{Symplectic triple systems}\label{se:STS}

The main objects of study in this section are defined as follows:

\begin{definition}\label{df:STS}
Let $T$ be a vector space endowed with a nonzero alternating
bilinear form $(.\vert.):T\times T\rightarrow k$, and a triple
product $T\times T\times T\rightarrow T$: $(x,y,z)\mapsto [xyz]$.
Then $\bigl(T,[...],(.\vert.)\bigr)$ is said to be a
\emph{symplectic triple system} (see \cite{YamAs}) if it satisfies
the following identitities:
\begin{subequations}\label{eq:STS}
\begin{align}
&[xyz]=[yxz]\label{eq:STSa}\\
&[xyz]-[xzy]=(x\vert z)y-(x\vert y)z+2(y\vert z)x\label{eq:STSb}\\
&[xy[uvw]]=[[xyu]vw]+[u[xyv]w]+[uv[xyw]]\label{eq:STSc}\\
&([xyu]\vert v)+(u\vert [xyv])=0\label{eq:STSd}
\end{align}
\end{subequations}
for any elements $x,y,z,u,v,w\in T$.
\end{definition}

Note that \eqref{eq:STSb} can be written as
\begin{equation}\label{eq:STSbb}
[xyz]-[xzy]=\psi_{x,y}(z)-\psi_{x,z}(y)
\end{equation}
with $\psi_{x,y}(z)=(x\vert z)y+(y\vert z)x$. (If $(.\vert.)$ is
nondegenerate, the maps $\psi_{x,y}$ span the symplectic Lie
algebra $\frsp(T)$.)

Also, \eqref{eq:STSc} is equivalent to $d_{x,y}=[xy.]$ being a
derivation of the triple system. Let $\inder(T)$ be the linear
span of $\{d_{x,y}:x,y\in T\}$, which is a Lie subalgebra of
$\End(T)$. Equation \eqref{eq:STSd} is equivalent to the condition
$d_{x,y}\in \frsp(T)$ for any $x,y\in T$.

As remarked in \cite[Section 4]{MSII}, equation \eqref{eq:STSd} is
a consequence of \eqref{eq:STSb} and \eqref{eq:STSc} unless the
characteristic of the ground field $k$ is $3$ and the dimension of
$T$ is $2$.

\medskip
As usual, a homomorphism of \STSs\ is a linear map
$\varphi:\T\rightarrow \bigl(T',[...]',(.\vert .)'\bigr)$
satisfying both $(\varphi(x)\vert \varphi(y))'=(x\vert y)$ and
$\varphi([xyz])=[\varphi(x)\varphi(y)\varphi(z)]'$ for any
$x,y,z\in T$. An ideal of a \STS\ is a subspace $I$ of $T$ such
that $[TTI]+[TIT]\subseteq I$, and the system is said to be simple
if $[TTT]\ne 0$ and it contains no proper ideal.

\begin{proposition}\label{pr:STSsimple}
Let $\T$ be a \STS\ and assume that the dimension of $T$ is $>2$
if the characteristic of the ground field $k$ is $3$. Then $\T$ is
simple if and only if the bilinear form $(.\vert.)$ is
nondegenerate.
\end{proposition}
\begin{proof} Let $T^\perp=\{x\in T: (x\vert T)=0\}$ be the kernel
of the bilinear map. Then by \eqref{eq:STSd} $[TTT^\perp]\subseteq
T^\perp$ and, by \eqref{eq:STSb}, $[TT^\perp T]\subseteq
[TTT^\perp] + T^\perp\subseteq T^\perp$. Hence $T^\perp$ is an
ideal of $T$. Therefore, if $T$ is simple, then $T^\perp =0$ and
$(.\vert .)$ is nondegenerate.

Conversely, assume that $(.\vert .)$ is nondegenerate and let $I$
be an ideal of $T$, then by \eqref{eq:STSb}
\[
(x\vert z)y-(x\vert y)z+2(y\vert z)x\in I
\]
for any $x,y,z\in T$ such that one of them is in $I$. With $y,z\in
T$ and $x\in I$, it follows that
\begin{equation}\label{eq:Iyz}
(x\vert z)y-(x\vert y)z\in I,
\end{equation}
while for $y\in I$ and $x,z\in T$ we get, after interchanging $x$
and $y$,
\begin{equation}\label{eq:Ixz}
(x\vert y)z+2(x\vert z)y\in I
\end{equation}
for any $y,z\in T$ and $x\in I$. If the characteristic of $k$ is
$\ne 3$, it follows from \eqref{eq:Iyz} and \eqref{eq:Ixz} that
$(I\vert T)T\subseteq I$, so that either $I\subseteq T^\perp =0$
or $I=T$. On the other hand, if $\charac k=3$ and $I\ne 0$, for
any $0\ne x\in I$, let $z\in T$ such that $(x\vert z)\ne 0$, then
\eqref{eq:Iyz} shows that $(kx)^\perp \subseteq I$. Therefore
either $I=T$ or $I$ has codimension one. But in the latter case,
if $\dim T>2$, then there are linearly independent elements
$x,y\in I$, and hence $T=(kx)^\perp + (ky)^\perp\subseteq I$.
Thus, either $I=T$ or $\dim T=2$, as required.
\end{proof}

\medskip

Proposition \ref{pr:STSsimple} fails for two dimensional \STSs\
over fields of characteristic $3$. Note that in this case, the
expression $(x\vert z)y-(x\vert y)z+2(y\vert z)x$ in the right
hand side of \eqref{eq:STSb} becomes $(x\vert z)y+(y\vert
x)z+(z\vert y)z$, which is identically zero (since it is
skew-symmetric in three variables on a two dimensional vector
space). Therefore, in this case, \eqref{eq:STSa} and
\eqref{eq:STSb} simply say that the triple product is symmetric.
The two dimensional \STSs\ in characteristic $3$ are easily
classified:

\begin{proposition}\label{pr:STS2} Let $\T$ be a two dimensional
\STS\ over a field $k$ of characteristic $3$. Then either:
\begin{romanenumerate}
\item There is a symplectic basis $\{a,b\}$ of $T$ (that is,
$(a\vert b)=1$) and a scalar $\alpha\in k$ such that $[aaa]=\alpha
b$ and all the other triple products of basis elements is $0$; or
\item there is a symplectic basis $\{a,b\}$ of $T$ and a scalar
$0\ne\epsilon\in k$ such that $[aaa]=0=[bbb]$ and
\begin{equation}\label{eq:STS2}
[aab]=[aba]=[baa]=\epsilon a,\quad [abb]=[bab]=[bba]=-\epsilon b.
\end{equation}
\end{romanenumerate}
Moreover, two \STSs\ in case (i) with scalars $\alpha$ and
$\alpha'$ are isomorphic if and only if there is a scalar
$0\ne\mu\in k$ with $\alpha=\mu^4\alpha'$, while the scalar
$\epsilon$ is an invariant of the isomorphism class of the \STS\
in case (ii).
\end{proposition}
\begin{proof}
The case of trivial product ($[TTT]=0)$ is covered in item (i).
Hence we assume that $[TTT]\ne 0$.

If there is an element $a\in T$ such that $[aaa]\ne 0$, take
$c=[aaa]$. Then by (\ref{eq:STS}a-c), for any $x,y\in T$,
$[xx[yyy]]=3[[xxy]yy]=0$, so that $[xxc]=0$ for any $x\in T$ and
hence $[TTc]=0$. In particular, this shows that $c\in T\setminus
ka$. Then there is a scalar $0\ne \alpha\in k$ such that $(a\vert
c)=\alpha$, and hence we get (i) with $b=\alpha^{-1}c$.

Therefore we may assume now that $[xxx]=0$ for any $x\in T$ but
$[TTT]\ne 0$. Then there is an element $a\in T$ such that
$[aaT]\ne 0$. If $d_{a,a}=[aa.]$ is not nilpotent, since $[aaa]=0$
there is an element $b\in T\setminus ka$ such that $[aab]=\lambda
b$ for some $0\ne \lambda\in k$. Thus, using \eqref{eq:STSc},
\[
\begin{split}
\lambda[abb]=[ab[aab]]&=2[[aab]ab]+[aa[abb]]\\
&=2\lambda[abb]+[aa[abb]],
\end{split}
\]
so that $[aa[abb]]=-\lambda[abb]$ and $[abb]=0$ as $0$ and
$\lambda$ are the only eigenvalues of $d_{a,a}$ ($\dim T=2$). In
this case, $0\ne \lambda (a\vert b)=(a\vert [aab])=-([aaa]\vert
b)=0$, a contradiction. Hence we conclude that $d_{a,a}$ is
nilpotent, so there is an element $c\in T\setminus ka$ such that
$[aac]=a$. From \eqref{eq:STS} it follows that
\[
[aa[acc]]=2[a[aac]c]=-[aac],
\]
so $[acc]+c\in\ker d_{a,a}=ka$, and there is a scalar $\mu\in k$
such that $[acc]=\mu a-c$. Now
\[
\begin{split}
\mu a-c&=[cca]=[cc[aac]]=2[[cca]ac]\\
 &=-[(\mu a -c)ac]=-\mu a+(\mu a-c)=-c.
\end{split}
\]
Hence $\mu =0$ and $[cca]=-c$. If $(a\vert c)=\alpha$, with
$b=\alpha^{-1}c$ we obtain case (ii) with $\epsilon=\alpha^{-1}$.

Finally, unless $[TTT]=0$, the ideal $kb=[TTT]$ in case (i) is
fixed under any automorphism, so any other symplectic basis with a
similar multiplication is of the form $\{a'=\mu a+\eta
b,b'=\mu^{-1}b\}$ for some $\mu,\eta\in k$ with $\mu\ne 0$. In
this basis, $[a'a'a']=\mu^4\alpha b'$, whence the condition for
isomorphism for these algebras. Also, in case (ii),
$d_{x,y}=-\epsilon\psi_{x,y}=-\epsilon\Bigl((x\vert .)y+(y\vert
.)x\Bigr)$ for any $x,y\in T$. Hence $\epsilon$ is an invariant of
the isomorphism class.
\end{proof}

\medskip

Note that the \STS\ in case (i) of Proposition \ref{pr:STS2} is
not simple, even though $(.\vert .)$ is nondegenerate.

\medskip

The simple \STS\ in case (ii) will be denoted by $T_{2,\epsilon}$.

\bigskip

Symplectic triple systems are strongly related to $\bZ_2$-graded
Lie algebras (or to Lie triple systems \cite{YamAs}). The precise
statement that will be most useful here is the following:

\begin{theorem}\label{th:STSLie}
Let $\bigl(T,[...],(.\vert .)\bigr)$ be a symplectic triple system
and let $\bigl(V,\langle.\vert.\rangle\bigr)$ be a two dimensional
vector space endowed with a nonzero alternating bilinear form.
Define the $\bZ_2$-graded algebra
$\frg=\frg(T)=\frgo\oplus\frguno$ with
\[
\begin{cases}
\frgo= \frsp(V)\oplus \inder(T)&\text{(direct sum of ideals),}\\
\frguno=V\otimes T\,,
\end{cases}
\]
and anticommutative multiplication given by:
\begin{itemize}
\item
$\frgo$ is a Lie subalgebra of $\frg$,
\item
$\frgo$ acts naturally on $\frguno$; that is
\[
[s,v\otimes x]=s(v)\otimes x,\qquad [d,v\otimes x]=v\otimes d(x),
\]
for any $s\in \frsp(V)$, $d\in\inder(T)$, $v\in V$, and $x\in T$.
\item
For any $u,v\in V$ and $x,y\in T$:
\begin{equation}\label{eq:oddbracket}
[u\otimes x,v\otimes y]=
  (x\vert y)\gamma_{u,v} +\langle u\vert v\rangle d_{x,y}
\end{equation}
where $\gamma_{u,v}=\langle u\vert .\rangle v+\langle v\vert
.\rangle u$ and $d_{x,y}=[xy.]$.
\end{itemize}
Then $\frg(T)$ is a Lie algebra. Moreover, $\frg(T)$ is simple if
and only if so is $\T$.

Conversely, given a $\bZ_2$-graded Lie algebra
$\frg=\frgo\oplus\frguno$ with
\[
\begin{cases}
\frgo=\frsp(V)\oplus \frs &\text{(direct sum of ideals),}\\
\frguno=V\otimes T&\text{(as a module for $\frgo$),}
\end{cases}
\]
where $T$ is a module for $\frs$, the $\frgo$-invariance of the
Lie bracket implies that equation \eqref{eq:oddbracket} is
satisfied for an alternating bilinear form $(.\vert .):T\times
T\rightarrow k$ and a symmetric bilinear map $d_{.,.}:T\times
T\rightarrow \frs$. Then, if $(.\vert .)$ is not $0$ and a triple
product on $T$ is defined by means of $[xyz]=d_{x,y}(z)$,
$\bigl(T,[...],(.\vert .)\bigr)$ is a symplectic triple system.
\end{theorem}
\begin{proof}
This is proved in \cite{MSII}, following the ideas in
\cite{YamAs}, with the exception of the assertion about the
simplicity of $\frg$. For this, first note that if $T$ is not
simple, then $T^\perp$ is a proper ideal of $T$ (see
\ref{pr:STSsimple} and its proof), and then $d_{T,T^\perp}$ is an
ideal of $\inder(T)$ and $\fra=d_{T,T^\perp}\oplus \bigl( V\otimes
T^\perp\bigr)$ is a proper ideal of $\frg$. Conversely, assume
that $T$ is simple. Since $\frguno$ is not the adjoint module for
$\frgo$, it is enough to prove that $\frg$ has no proper
homogeneous ideal. Let $\fra=\fra\subo\oplus \fra\subuno$ be an
homogeneous ideal of $\frg$. If $\fra\subuno=0$, then $\fra$ is an
ideal of $\frgo$ and $[\fra,\frguno]=0$. But $\frgo$ acts
faithfully on $\frguno$, hence $\fra=0$. Otherwise,
$\fra\subuno\ne 0$ so, by $\frsp(V)$-invariance,
$\fra\subuno=V\otimes I$ for some subspace $I$ of $T$, which is
closed under the action of $\inder(T)$. Thus, $[TTI]\subseteq I$.
By $\frsp(V)$-invariance too, $\fra\subo=\bigl(\fra\subo\cap
\frsp(V)\bigr) \oplus \bigl(\fra\subo\cap\inder(T)\bigr)$. Let
$\{v,w\}$ be a symplectic basis of $V$, then for any $x\in T$ and
$y\in I$, $[x\otimes v,y\otimes w]=d_{x,y}+(x\vert
y)\gamma_{v,w}\in\fra\subo$. Hence $d_{I,T}\subseteq \fra$ and,
hence, $[d_{I,T},\frguno]=[d_{I,T},V\otimes I]=V\otimes
[ITT]\subseteq \fra$. Thus $[ITT]\subseteq I$ and $I$ is an ideal
of $T$. But then $I=T$, $\fra\subuno=\frguno$ and
$\frgo=[\frguno,\frguno]\subseteq \fra$, so that $\fra=\frg$.
\end{proof}

\medskip

The argument in the proof above is similar to the one in
\cite[Theorem 2.2]{EKO}.

\medskip

As a noteworthy example, consider the simple \STS\
$T_{2,\epsilon}$ that appears in Proposition \ref{pr:STS2}, where
the characteristic of the ground field $k$ is $3$. Here, as
remarked in the proof of \ref{pr:STS2},
$d_{x,y}=-\epsilon\psi_{x,y}$ for any $x,y\in T_{2,\epsilon}$, so
that $\inder(T)=\frsp(T)$. Therefore, the simple Lie algebra
$\frg(T_{2,\epsilon})$ is given by:
\begin{equation}\label{eq:gt2e}
\frg(T_{2,\epsilon})=\Bigl(\frsp(V_1)\oplus\frsp(V_2)\Bigr)\oplus
\Bigl(V_1\otimes V_2\Bigr)
\end{equation}
where $V_1$ and $V_2$ are two dimensional vector spaces endowed
with nonzero alternating bilinear forms $\langle .\vert
.\rangle_i$ ($i=1,2$) and the Lie bracket is determined by the Lie
multiplication in $\frsp(V_1)\oplus \frsp(V_2)$, the natural
action of $\frsp(V_1)\oplus\frsp(V_2)$ on $V_1\otimes V_2$ and by
\[
[a\otimes x,b\otimes y]=\langle x\vert
y\rangle_2\gamma_{a,b}-\epsilon\langle a\vert
b\rangle_1\gamma_{x,y},
\]
for any $a,b\in V_1$ and $x,y\in V_2$, where $\gamma_{a,b}=\langle
a\vert .\rangle_1 b+ \langle b\vert .\rangle_1a$ and
$\gamma_{x,y}=\langle x\vert .\rangle_2 y+ \langle y\vert
.\rangle_2x$.

\smallskip

There is a great similarity of the above defined Lie algebra and
the exceptional simple classical Lie superalgebras
$\Gamma(\sigma_1,\sigma_2,\sigma_3)$ in \cite[pp.
17-18]{Scheunert}, where three two dimensional space are used.
These Lie superalgebras constitute the family $D(2,1;a)$ in Kac's
classification \cite{Kac}.

\smallskip

Kostrikin defined in \cite{Kostrikin} a parametric family
$L(\epsilon)$ ($\epsilon\ne 0$) of ten dimensional simple Lie
algebras over fields of characteristic $3$, which appear as
subalgebras of the contact algebra $K_3$. Here these simple Lie
algebras have a natural interpretation:

\begin{proposition}\label{pr:Lepsilon}
Let $k$ be a field of characteristic $3$. Then, for any $0\ne
\epsilon\in k$, the Lie algebra $\frg(T_{2,\epsilon})$ and
$L(\epsilon)$ are isomorphic.
\end{proposition}
\begin{proof} This has been checked in \cite{BrownCont}, although
in terms of Freudenthal triple systems instead of \STSs\ (see
Theorem \ref{th:STSFTS} below). Actually, the Freudenthal triple
system that appears in the second paragraph of
\cite[p.~29]{BrownCont} with parameter $c$ corresponds through
\ref{th:STSFTS} to the \STS\ $T_{2,\epsilon}$ with
$\epsilon=-(c+1)$, Brown shows that this system corresponds to the
Kostrikin algebra $L\left(-\frac{1}{c+1}\right)$, which in turn is
isomorphic to $L(-(c+1))=L(\epsilon)$ (\cite{Kostrikin}).
\end{proof}

\medskip

\begin{remark}\label{re:Lepsilon}
Therefore, the definition of $\frg(T_{2,\epsilon})$ in
\eqref{eq:gt2e} gives a nice description of the simple Lie algebra
$L(\epsilon)$. Also, take symplectic bases $\{v_i,w_i\}$ of $V_i$
($i=1,2$) and the elements
\begin{align*}
e_1&=\tfrac{1}{2}\gamma_{v_2,v_2},& f_1&=-\gamma_{w_2,w_2},&
h_1&=[e_1,f_1]=-\gamma_{v_2,w_2},\\
e_2&=v_1\otimes w_2,& f_2&=-w_1\otimes v_2,&
h_2&=[e_2,f_2]=\gamma_{v_1,w_1}+\epsilon\gamma_{v_2,w_2},
\end{align*}
which satisfy that
\[
[h_1,e_1]=2e_1,\ [h_1,e_2]=-e_2,\ [h_2,e_1]=-2\epsilon e_1,\
[h_2,e_2]=(\epsilon -1)e_2,
\]
and this shows that $\frg(T_{2,\epsilon})$ is the contragredient
Lie algebra (see \cite[\S 3]{VK}) associated to the Cartan matrix
$A=\bigl(\begin{smallmatrix} 2&-1\\ -2\epsilon&\epsilon
-1\end{smallmatrix}\bigr)=\bigl(\begin{smallmatrix} 2&-1\\
\epsilon&\epsilon -1\end{smallmatrix}\bigr)$.

For $\epsilon=1$ ($L(1)$ is shown in \cite{Kostrikin} to present
some specific features), after scaling the last row, this is the
matrix
$\textbf{C}_{2,\infty} = \bigl(\begin{smallmatrix} 2&-1\\
-1&0\end{smallmatrix}\bigr)$, while for $\epsilon \ne 1$, after
scaling the last row, we get the Cartan matrix
$\textbf{C}_{2,\frac{\epsilon}{1-\epsilon}}=\bigl(\begin{smallmatrix} 2&-1\\
\frac{\epsilon}{1-\epsilon}&2\end{smallmatrix}\bigr)$. In
particular, for $\epsilon=-1$, we get the Cartan matrix associated
to the simple classical Lie algebra of type $C_2$ and, indeed,
$\frg(T_{2,-1})$ is naturally isomorphic to $\frsp(V_1\perp V_2)$.
\qed
\end{remark}

\bigskip

For later use, let us define an eight dimensional module for the
simple Lie algebra $L(1)=\frg(T_{2,1})$ ($\charac k=3$), which
will provide an important example of a simple \STS. With the
previous notation, consider the $\bZ_2$-graded vector space
$W=W\subo\oplus W\subuno$, with
\[
\begin{cases}
W\subo =V_1\otimes\frsp(V_2),&\\
W\subuno=V_2,&
\end{cases}
\]
and give $W$ the structure of a graded module for the
$\bZ_2$-graded Lie algebra $\frg(T_{2,1})$ with the natural action
of $\frsp(V_1)\oplus\frsp(V_2)$ (that is, $\rho_{s_1}(a_1\otimes
s_2)=s_1(a_1)\otimes s_2$, $\rho_{s_2}(a_1\otimes t_2)=a_1\otimes
[s_2,t_2]$, $\rho_{s_1}(a_2)=0$ and $\rho_{s_2}(a_2)=s_2(a_2)$,
for any $a_1\in V_1$, $a_2\in V_2$, $s_1\in \frsp(V_1)$ and
$s_2,t_2\in \frsp(V_2)$). Finally, for any $a_1,b_1\in V_1$,
$a_2,b_2\in V_2$ and $s_2\in \frsp(V_2)$, define the action of
$a_1\otimes a_2\in \frg(T_{2,1})\subuno$ on $W$ by means of:
\[
\begin{split}
&\rho_{a_1\otimes a_2}(b_1\otimes s_2)=-\langle a_1\vert
b_1\rangle_1s_2(a_2)\in W\subuno,\\
&\rho_{a_1\otimes a_2}(b_2)=-a_1\otimes \gamma_{a_2,b_2}\in
W\subo.
\end{split}
\]
Using that
$[s_2,\gamma_{a_2,b_2}]=\gamma_{s_2(a_2),b_2}+\gamma_{a_2,s_2(b_2)}$
and that, since the characteristic is $3$, $(a_2\vert
b_2)s_2=\gamma_{s_2(a_2),b_2}-\gamma_{a_2,s_2(b_2)}$ for any
$a_2,b_2\in V_2$ and $s_2\in \frsp(V_2)$, it follows easily that
$\rho:\frg(T_{2,1})\rightarrow \frgl(W)$ is indeed a
representation.

\begin{proposition}\label{pr:STS8}
Let $k$ be a field of characteristic $3$ and let $\rho:
\frg(T_{2,1})\rightarrow \frgl(W)$ be the representation defined
above. Then $T=W$, endowed with the alternating bilinear form
given by
\[
\left\{%
\begin{aligned}
&(W\subo\vert W\subuno)=0,\\
&(a_1\otimes s_2\vert b_1\otimes t_2)=\langle a\vert
b\rangle_1\trace(s_2t_2),\\
&(a_2\vert b_2)=\langle a_2\vert b_2\rangle_2,
\end{aligned}\right.
\]
and with the triple product $[ABC]=d_{A,B}(C)$, where $d:W\times
W\rightarrow \frgl(W)$, $(A,B)\mapsto d_{A,B}$ is the symmetric
bilinear map determined by
\[
\left\{%
\begin{aligned}
& d_{a_1\otimes s_2,b_1\otimes
t_2}=\rho_{\trace(s_2t_2)\gamma_{a_1,b_1}-\langle a_1\vert
b_1\rangle_1 [s_2,t_2]},\\
 &d_{a_1\otimes s_2,a_2}=\rho_{a_1\otimes s_2(a_2)},\\
 &d_{a_2,b_2}=-\rho_{\gamma_{a_2,b_2}},
\end{aligned}\right.
\]
for any $a_1,b_1\in V_1$, $a_2,b_2\in V_2$ and $s_2,t_2\in
\frsp(V_2)$, is a \STS.
\end{proposition}
\begin{proof}
This is a straightforward computation, using that for any
$s_2,t_2,r_2\in \frsp(V_2)$,
\[
\begin{split}
s_2t_2&=\frac{1}{2}\trace(s_2t_2)+\frac{1}{2}[s_2,t_2]
   =-\trace(s_2t_2)-[s_2,t_2]\in\frgl(V_2),\\
[[s_2,t_2],r_2]&=s_2(t_2r_2+r_2t_2)+(t_2r_2+r_2t_2)s_2\\
&\qquad\qquad \qquad\qquad -(s_2r_2+r_2s_2)t_2
 -t_2(r_2s_2+s_2r_2)\\
 &=2\trace(t_2r_2)s_2-2\trace(s_2r_2)t_2\\
 &=\trace(s_2r_2)t_2-\trace(t_2r_2)s_2. \qedhere
\end{split}
\]
\end{proof}

\bigskip

The \STSs\ are strongly related to other triple systems too, in
particular to Freudenthal triple systems and to Faulkner ternary
algebras \cite{Faulkner} (or balanced symplectic algebras
\cite{FaulknerFerrar}).

Let $T$ be a vector space endowed with a nondegenerate alternating
bilinear form $(.\vert .):T\times T\rightarrow k$, and a triple
product $T\times T\times T\rightarrow T$, $(x,y,z)\mapsto xyz$.
Then $\bigl(T,xyz,(.\vert .)\bigr)$ is said to be a
\emph{Freudenthal triple system} (see
\cite{Meyberg,Ferrar,BrownFTS}) if it satisfies:
\begin{subequations}\label{eq:FT}
\begin{align}
&\text{$xyz$ is symmetric in its arguments,}\label{eq:FTa}\\
&\text{$( x\vert yzt)$ is symmetric in its
arguments,}\label{eq:FTb}\\
&(xyy)xz+(yxx)yz+(xyy\vert z)x+(yxx\vert z)y\nonumber\\
&\hskip 2in +(x\vert z)xyy+(y\vert z)yxx=0, \label{eq:FTc}
\end{align}
\end{subequations}
for any $x,y,z,t\in T$.

The next result relates the \STSs\ and the Freudenthal triple
systems, its proof may be found in \cite[Theorem 4.7]{MSII}.

\begin{theorem}\label{th:STSFTS}
Let $(.\vert .)$ be a nondegenerate alternating bilinear form on
the vector space $T$ and let $xyz$ and $[xyz]$ be two triple
products on $T$ related by $xyz=[xyz]-\psi_{x,y}(z)$ for any
$x,y,z\in T$. Then $\T$ is a \STS\ if and only if either $xyz=0$
for any $x,y,z\in T$, or  $\bigl(T,xyz,(.\vert .)\bigr)$ is a
Freudenthal triple system.
\end{theorem}

\medskip

Also, let $T$ be a vector space endowed with an alternating
bilinear form $(.\vert .)$ and a triple product $\langle
xyz\rangle$ satisfying (\cite{Faulkner} or \cite[Section
3]{FaulknerFerrar}):
\begin{subequations}\label{eq:Faulkner}
\begin{align}
&\langle xyz\rangle=\langle yxz\rangle +(x\vert y)z\label{eq:Faulknera}\\
&\langle xyz\rangle =\langle xzy\rangle +(y\vert z)x\label{eq:Faulknerb}\\
&\langle\langle xyz\rangle vw\rangle=\langle\langle xvw\rangle
yz\rangle +\langle x\langle yvw\rangle z\rangle +\langle xy\langle
zwv\rangle\rangle \label{eq:Faulknerc}
\end{align}
\end{subequations}
for any $x,y,z,v,w\in T$. Then $\bigl(T,\langle ...\rangle,(.\vert
.)\bigr)$ is called a \emph{Faulkner ternary algebra}. These
triple systems have been called \emph{balanced symplectic
algebras} in \cite{FaulknerFerrar}.

\begin{theorem}\label{th:STSFaulkner}
Let $(.\vert .)$ be an alternating bilinear form on the vector
space $T$ and let $[xyz]$ and $\langle xyz\rangle$ be two triple
products related by $[xyz]=-2\langle zxy\rangle +(x\vert y)z$ for
any $x,y,z\in T$. Then $\T$ satisfies (\ref{eq:STS}a-d) if and
only if $\bigl(T,\langle ...\rangle,(.\vert .)\bigr)$ is a
Faulkner ternary algebra. In particular, the \STSs\ are in
bijection with the Faulkner ternary algebras with nonzero
alternating bilinear form.
\end{theorem}
\begin{proof}
For any $x,y,z\in T$,
\[
\begin{split}
\langle xyz\rangle-\langle yxz\rangle -(x\vert y)z
 &=-\frac{1}{2}\Bigl([yzx]-[xzy]-(z\vert x)y-(z\vert
 y)x\Bigr)-(x\vert y)z\\
 &=-\frac{1}{2}\Bigl([yzx]-[xzy]-\psi_{z,y}(x)+\psi_{x,z}(y)\Bigr),\\
\langle xyz\rangle-\langle xzy\rangle-(y\vert z)x
 &=-\frac{1}{2}\Bigl([yzx]-[zyx]-(y\vert z)x+(z\vert
 y)x\Bigr)-(y\vert z)x\\
 &=-\frac{1}{2}\Bigl([yzx]-[zyx]\Bigr),
\end{split}
\]
so that \eqref{eq:STSa} and \eqref{eq:STSb} are equivalent to
\eqref{eq:Faulknera} and \eqref{eq:Faulknerb}.

On the other hand, if \eqref{eq:Faulknerb} and
\eqref{eq:Faulknerc} are satisfied,  \cite[(3.1)]{FaulknerFerrar}
shows that
\begin{equation}\label{eq:Faulknerd}
(\langle xzw\rangle \vert y)+(x\vert \langle ywz\rangle )=0
\end{equation}
for any $x,y,z,t\in T$, but
\[
\begin{split}
-2\Bigl( (\langle xzw\rangle \vert y)+&(x\vert \langle
ywz\rangle)\Bigr)\\
 &=([zwx]\vert y)-(z\vert w)(x\vert y)+(x\vert [wzx])-(w\vert
 z)(x\vert y)\\
 &=([zwx]\vert y)+(x\vert [wzy]).
\end{split}
\]
Hence \eqref{eq:STSa} and \eqref{eq:STSd} are equivalent to
\eqref{eq:Faulknera} and \eqref{eq:Faulknerd}. Finally, assuming
\eqref{eq:STSa} and \eqref{eq:STSb}
\[
\begin{split}
4\Bigl( \langle\langle xyz\rangle vw\rangle &-\langle\langle
xvw\rangle
  yz\rangle -\langle x\langle yvw\rangle z\rangle -\langle xy\langle
   zwv\rangle\rangle\Bigr)\\
&= -2\Bigr( [vw\langle xyz\rangle]-(v\vert w)\langle xyz\rangle
-[yz\langle xvw\rangle]+(y\vert z)\langle xvw\rangle\\
&\qquad -[\langle yvw\rangle zx]+(\langle yvw\rangle\vert z)x-
   [y\langle zwx\rangle x]+(y\vert\langle zwv\rangle)x\Bigr)\\
&=[vw[yzx]]-(y\vert z)[vwx]-(v\vert w)\bigl([yzx]
 -(y\vert z)x\bigr)\\
 &\qquad -[yz[vwx]]+(v\vert w)[yzx]+(y\vert z)[vwx]-(y\vert
 z)(v\vert w)x\\
 &\qquad -[[vwy]zx]+(v\vert w)[yzx]+([vwy]\vert z)x-(v\vert
 w)(y\vert z)x\\
 &\qquad -[y[wvz]x]+(w\vert v)[yzx]+(y\vert[wvz])x-(y\vert
 z)(w\vert v)x\\
&=\Bigl([vw[yzx]]-[y[vwz]x]-[[vwy]zx]-[yz[vwx]]\Bigr)\\
  &\qquad\qquad\qquad+\Bigl(([vwy]\vert z)+(y\vert [vwz])\Bigr)x
\end{split}
\]
which shows that, if \eqref{eq:STSd} is satisfied, then
\eqref{eq:STSc} and \eqref{eq:Faulknerd} are equivalent.
\end{proof}

\begin{remark}
See \cite{BrownFTS} for a proof, in characteristic $3$, of the
relationship between Freudenthal triple systems and Faulkner
ternary algebras. Also, \cite[Lemma 3.2 and Corollary
1]{FaulknerFerrar} deals with this relationship in characteristic
$\ne 2,3$.
\end{remark}

\medskip

Now, Theorem \ref{th:STSFaulkner}, together with the
classification of the simple Faulkner ternary algebras with
nonzero alternating bilinear form in \cite[Theorem
4.1]{FaulknerFerrar} over fields of characteristic $\ne 2,3$
(which is based on the classification in \cite{Meyberg} of the
Freudenthal triple systems) immediately yields the classification
of the simple \STSs\ over algebraically closed fields. Here we
follow the notations in \cite{McCrimmon} concerning Jordan
algebras.

\begin{theorem}\label{th:charne23}
Let $k$ be an algebraically closed field of characteristic $\ne
2,3$ and let $\T$ be a simple symplectic triple system. Then
either:
\begin{romanenumerate}
\item $[xyz]=\psi_{x,y}(z)$ for any $x,y,z\in T$; or
\item there exists a Jordan algebras $J$ such that either $J=0$,
or $J=\calJord(q,e)$ is the Jordan algebra of a nondegenerate
quadratic form $q$ with basepoint $e$ (\cite[II.3.3]{McCrimmon}),
with trace form $t(a,b)$, where we define $a\times b=0$ for any
$a,b\in J$, or a Jordan algebra $J=\calJord(n,c)$ of a
nondegenerate cubic form $n$ with basepoint $c$
(\cite[II.4.3]{McCrimmon}), trace form $t(a,b)$ and cross product
$a\times b$, such that, up to isomorphism,
\begin{equation}\label{eq:FJJF}
T=\left\{\begin{pmatrix} \alpha &a\\ b&\beta\end{pmatrix} :
\alpha,\beta\in k,\ a,b\in J\right\},
\end{equation}
and for $x_i=\bigl(\begin{smallmatrix} \alpha_i&a_i\\
b_i&\beta_i\end{smallmatrix}\bigr)$, $i=1,2,3$:
\begin{equation}\label{eq:FJJFbis}
\left\{%
\begin{aligned}
&(x_1\vert x_2)=
  \alpha_1\beta_2-\alpha_2\beta_1-t(a_1,b_2)+t(b_1,a_2),\\
&[x_1x_2x_3]=\begin{pmatrix} \gamma&c\\ d&\delta
    \end{pmatrix}\quad\text{with}\\
&\quad \gamma=\Bigl(-3(\alpha_1\beta_2+\alpha_2\beta_1)
   +t(a_1,b_2)+t(a_2,b_1)\Bigr)\alpha_3\\
&\quad\quad +2\Bigl(\alpha_1t(b_2,a_3)
   +\alpha_2t(b_1,a_3)-t(a_1\times a_2,a_3)\Bigr)\\
&\quad c=\Bigl(-(\alpha_1\beta_2+\alpha_2\beta_1)
   +t(a_1,b_2)+t(a_2,b_1)\Bigr)a_3\\
&\quad\quad +2\Bigl(\bigl(t(b_2,a_3)-\beta_2\alpha_3\bigr)a_1+
    \bigl(t(b_1,a_3)-\beta_1\alpha_3\bigr)a_2\Bigr)\\
&\quad\quad +2\Bigl( \alpha_1(b_2\times b_3)
   +\alpha_2(b_1\times b_3)+\alpha_3(b_1\times b_2)\Bigr)\\
&\quad\quad -2\Bigl( (a_1\times a_2)\times b_3 +
  (a_1\times a_3)\times b_2 +(a_2\times a_3)\times b_1\Bigr)\\
&\quad \text{$\delta=-\gamma^\sigma$, $d=-c^\sigma$, where
$\sigma=(\alpha\beta)(ab)$ (that is, $\gamma^{\sigma}$ and
$c^\sigma$}\\
&\quad\quad\text{are obtained from $\gamma$ and $c$ by
interchanging $\alpha$ and $\beta$ and}\\
&\quad\quad\text{also $a$ and $b$ throughout).} \qed
\end{aligned}
\right.
\end{equation}
\end{romanenumerate}
\end{theorem}

\bigskip

\begin{remark}
The case $J=0$ above is missing in \cite[Theorem
4.1]{FaulknerFerrar}, it corresponds to the case in which
$\mathcal{M}^+=0=\mathcal{M}^-$ in the arguments of \cite[Section
4]{FaulknerFerrar}.
\end{remark}

\begin{remark}\label{re:Jordan3}
Over an algebraically closed ground field $k$ of characteristic
$\ne 2,3$, the Jordan algebras $\calJord(n,c)$ of a nondegenerate
cubic form with basepoint are (see \cite{JacobsonJordan,
McCrimmon}), up to isomorphism, either the ground field $k$ with
$n(\alpha)=\alpha^3$, so that $t(\alpha,\beta)=3\alpha\beta$ for
any $\alpha,\beta\in k$, or the cartesian product $k\times
\calJord(q,e)$, for a nondegenerate quadratic form $q$, or the
Jordan algebra $J=H_3(C)$ of $3\times 3$-hermitian matrices over
$C$, where $C$ is either $k$, $k\times k$, $\Mat_2(k)$ or the
algebra of split octonions. In all cases, $c$  is the identity
element $1$.
\end{remark}

\medskip

A natural question now is, given the \STSs\ in Theorem
\ref{th:charne23}, what do the simple Lie algebras $\frg(T)$ look
like? In order to answer this question, and for future use, we
will first consider different constructions of \STSs.

\begin{examples}\label{ex:examples}
Let $\bigl(V,\langle.\vert.\rangle\bigr)$ be a two dimensional
vector space endowed with a nonzero alternating bilinear form.

\medskip

\noindent\emph{Symplectic case:}\quad Let $\T$ be a \STS\ in item
(i) of Theorem \ref{th:charne23}, then the Lie algebra
$\frg(T)=\frgo\oplus\frguno$ defined in \ref{th:STSLie} satisfies
that $\frgo=\frsp(V)\oplus\frsp(T)$ and $\frguno=V\otimes T$.
Moreover, $\frgo$ embeds naturally in the symplectic Lie algebra
$\frsp(V\perp T)$ of the orthogonal direct sum $V\perp T$, and so
does $\frguno=V\otimes T$ by means of $u\otimes x\rightarrow
\Psi_{u,x}: v\mapsto \langle u\vert v\rangle x,\ y\mapsto (x\vert
y)u$, for any $u,v\in V$ and $x,y\in T$. This gives an isomorphism
$\frg(T)\simeq \frsp(V\perp T)$. \qed

\medskip

\noindent\emph{Special case:}\quad Let $W$ be a nonzero vector
space over a ground field $k$ and let $W^*$ be its dual vector
space. Consider the direct sum $T=W\oplus W^*$, endowed with the
triple product determined by $[WWT]=0=[W^*W^*T]$ and
\begin{subequations}\label{eq:xfyg}
\begin{align}
[xfy]&=f(x)y+2f(y)x,\label{eq:xfy}\\
[xfg]&=-f(x)g-2g(x)f,\label{eq:xfg}
\end{align}
\end{subequations}
for any $x,y\in W$ and $f,g\in W$, and of the alternating bilinear
form such that $(W\vert W)=0=(W^*\vert W^*)$ and $(f\vert
x)=-(x\vert f)=f(x)$ for any $x\in W$ and $f\in W^*$. Then $\T$ is
easily checked to be a \STS, where $W$ and $W^*$ are invariant
under $\inder(T)$. Also note that for any $x\in W$ and $f\in W^*$,
the element $d_{x,f}\vert_W\in\frgl(W)$ (determined by
\eqref{eq:xfy}) has trace equal to $(2+\dim W)f(x)$, which is $0$
if the characteristic of $k$ divides $2+\dim W$.

Moreover, let $\frg(T)=\bigl(\frsp(V)\oplus\inder(T)\bigr)\oplus
(V\otimes T)$ be the associated $\bZ_2$-graded Lie algebra
(Theorem \ref{th:STSLie}), and identify the Lie algebra
$\frgl(V\oplus W)$ of endomorphisms of $V\oplus W$ with the set of
$2$ by $2$
matrices $\bigl(\begin{smallmatrix} A&B\\
C&D\end{smallmatrix}\bigr)$, where $A\in \frgl(V)$,
$D\in\frgl(W)$, $B\in\Hom_k(W,V)$ and $C\in\Hom_k(V,W)$. Then  if
$\frz$ denotes the center of the special linear Lie algebra
$\frsl(V\oplus W)$ ($\frz=0$ if $\charac k$ dos not divide
$\dim(V\oplus W)=2+\dim W$, and $\frz=k\, id$ otherwise) and
$\frpsl(V\oplus W)=\frsl(V\oplus W)/\frz$ denotes the
corresponding projective special linear Lie algebra, a
straightforward computation shows that the linear map $
\Phi:\frg(T)\rightarrow \frpsl(V\oplus W)$ such that
\[
\begin{split}
\Phi(s)&=\begin{pmatrix} s&0\\ 0&0\end{pmatrix} +\frz\,,\\
\Phi(d)&=\begin{cases} \begin{pmatrix} 0&0\\
0&d\end{pmatrix}+\frz\qquad\text{if $\charac k$ divides $2+\dim W$,}\\
\noalign{\smallskip}
\begin{pmatrix} 0&0\\ 0&d\end{pmatrix} -\frac{\trace(d)}{2+\dim
W}\begin{pmatrix} id&0\\ 0&id\end{pmatrix}\qquad\text{otherwise},
\end{cases}\\
\Phi(a\otimes x)&=\begin{pmatrix} 0&0\\ \langle a\vert .\rangle
x&0\end{pmatrix} +\frz\,,\\
\Phi(a\otimes f)&=\begin{pmatrix} 0&2f(.)a\\ 0&0\end{pmatrix}
+\frz\,,
\end{split}
\]
for any $s\in\frsp(V)$, $d\in \inder(T)$, $a\in V$, $x\in W$ and
$f\in W^*$, is a Lie algebra isomorphism. \qed

\medskip

\noindent\emph{Orthogonal case:}\quad Consider now a vector space
$W$ of dimension $\geq 3$ over our ground field $k$, endowed with
a nondegenerate symmetric bilinear form $q_W(.,.)$. On the tensor
product $V\otimes V$ there is another nondegenerate symmetric
bilinear form determined by $q_{V\otimes V}(a\otimes u,b\otimes
v)=\langle a\vert b\rangle\langle u\vert v\rangle$.

Recall that if $(U,q)$ is a vector space endowed with a
nondegenerate bilinear form, then the orthogonal Lie algebra
\[
\frso(U,q)=\left\{ f\in \frgl(U): q(f(x),y)+q(x,f(y))=0,\ \forall
x,y\in U\right\}
\]
is spanned by the linear operators $\sigma_{x,y}: z\mapsto
q(x,z)y-q(y,z)x$, for $x,y\in U$, which satisfy that
$[\sigma_{x,y},\sigma_{z,t}]=\sigma_{\sigma_{x,y}(z),t}+\sigma_{z,\sigma_{x,y}(t)}$
for any $x,y,z,t\in U$. Moreover, if $U$ is the orthogonal direct
sum $U=U_1\perp U_2$ of two subspaces, then $\frso(U_i,q_i)$
$(q_i=q\vert_{U_i}$) is naturally embedded in $\frso(U,q)$,
$i=1,2$, and $\frso(U,q)=\frso(U_1,q_1)\oplus\frso(U_2,q_2)\oplus
\sigma_{U_1,U_2}$. Besides, $\sigma_{U_1,U_2}$ is linearly
isomorphic to $U_1\otimes U_2$ and
\[
\begin{split}
[\sigma_{x_1,x_2},\sigma_{y_1,y_2}]&=
 \sigma_{\sigma_{x_1,x_2}(y_1),y_2}+\sigma_{y_1,\sigma_{x_1,x_2}(y_2)}\\
 &=q_2(x_2,y_2)\sigma_{x_1,y_1}+q_1(x_1,y_1)\sigma_{x_2,y_2}\in
 \frso(U_1,q_1)\oplus\frso(U_2,q_2)
\end{split}
\]
for any $x_1,y_1\in U_1$ and $x_2,y_2\in U_2$.

The vector space $(V\otimes V)\oplus W$ is endowed with the
nondegenerate symmetric bilinear form which is the orthogonal sum
of $q_{V\otimes V}$ and $q_W$, and hence its orthogonal Lie
algebra is
\[
\frso\bigl((V\otimes V)\oplus W\bigr)=\frso(V\otimes V,q_{V\otimes
V})\oplus \frso(W,q_W)\oplus \sigma_{V\otimes V,W}.
\]
But $\frso(V\otimes V,q_{V\otimes V})=\frsp(V)\oplus \frsp(V)$,
where the first (respectively, second) copy of $\frsp(V)$ acts on
the first (resp., second) slot of $V\otimes V$. Therefore, with
the natural identifications, one has
\[
\frso\bigl((V\otimes V)\oplus
W\bigr)=\Bigl(\frsp(V)\oplus\frsp(V)\oplus
\frso(W,q_W)\Bigr)\oplus (V\otimes V\otimes W),
\]
and
\begin{equation}\label{eq:orthogonal}
[a\otimes u\otimes x,b\otimes v\otimes y]=q_W(x,y)\sigma_{a\otimes
u,b\otimes v}+q_{V\otimes V}(a\otimes u,b\otimes v)\sigma_{x,y},
\end{equation}
for any $a,b,u,v\in V$ and $x,y\in W$. But $q_{V\otimes
V}(a\otimes u,b\otimes v)=\langle a\vert b\rangle \langle u\vert
v\rangle$ and
\[
\begin{split}
\sigma_{a\otimes u,b\otimes v}(c\otimes w)&=
 \langle a\vert c\rangle\langle u\vert w\rangle b\otimes v-
 \langle b\vert c\rangle\langle v\vert w\rangle a\otimes u\\
 &=\frac{1}{2}\langle u\vert v\rangle \gamma_{a,b}(c)\otimes w
 +\frac{1}{2}\langle a\vert b\rangle c\otimes \gamma_{u,v}(w),
\end{split}
\]
for any $a,b,c,u,v,w\in V$. (To check the last equality, it is
enough by Zariski density to assume that $\{a,b\}$ and $\{u,v\}$
are bases of $V$, and then even that $\langle a\vert b\rangle
=1=\langle u\vert v\rangle$. Now one can just express $c$
(respectively $w$) as a linear combination of $a$ and $b$ (resp.,
of $u$ and $v$) and expand.)

Hence equation \eqref{eq:orthogonal} becomes
\[
[a\otimes u\otimes x,b\otimes v\otimes y]
 =\frac{1}{2}\langle u\vert v\rangle q_W(x,y)\gamma_{a,b}+
 \frac{1}{2}\langle a\vert b\rangle q_W(x,y)\gamma_{u,v}+
 \langle a\vert b\rangle \langle u\vert v\rangle \sigma_{x,y},
\]
for any $a,b,u,v\in V$ and $x,y\in W$. According to Theorem
\ref{th:STSLie}, $T=V\otimes W$ is a  \STS, endowed with the
nondegenerate alternating form given by $\bigl(u\otimes x\vert
v\otimes y\bigr)=\frac{1}{2}\langle u\vert v\rangle q_W(x,y)$, and
with the triple product given by
\begin{equation}\label{eq:STSorthogonal}
[(u\otimes x)(v\otimes y)(w\otimes z)]=
 \frac{1}{2}q_W(x,y)\gamma_{u,v}(w)\otimes z +
 \langle u\vert v\rangle \sigma_{x,y}(z),
\end{equation}
for any $u,v,w\in V$ and $x,y,z\in T$. Besides, by Proposition
\ref{pr:STSsimple}, this \STS\ is simple. Moreover, the Lie
algebra $\frg(T)$ is isomorphic to the orthogonal Lie algebra
$\frso\bigl((V\otimes V)\oplus W,q\bigr)$. \qed

\medskip

\noindent\emph{G$_2$-case:}\quad Assume here that the
characteristic of the ground field $k$ is $\ne 3$. Let $k[x,y]$ be
the polynomial algebra in two variables $x$ and $y$, and identify
$\frsl_2(k)\,(\simeq \frsp(V))$ with the subalgebra of the Lie
algebra of derivations of $k[x,y]$ spanned by
$\left\{x\frac{\partial\ }{\partial x}-y\frac{\partial\ }{\partial
y},x\frac{\partial\ }{\partial y},y\frac{\partial\ }{\partial
x}\right\}$. Let $V_n=k_n[x,y]$ be the linear space of the degree
$n$ homogeneous polynomials. Then $V_n$ is invariant under the
action of $\frsl_2(k)$ and if the characteristic of $k$ is either
$0$ or $>n$, then it is irreducible. In particular, since $\charac
k\ne 3$, $V_3$ is an irreducible module of dimension $4$ for
$\frsl_2(k)$.

For any $f\in V_n$ and $g\in V_m$, the \emph{transvection}
$(f,g)_q$ is defined, assuming no conflict with the characteristic
of $k$, by \cite{Dixmier}:
\[
(f,g)_q=\begin{cases} 0\quad\text{if $q>\min(n,m)$,}\\
  \tfrac{(n-q)!}{n!}\tfrac{(m-q)!}{m!}\sum_{i=0}^q\left(
  (-1)^i\frac{\partial^qf}{\partial x^{q-i}\partial y^i}
  \frac{\partial^qg}{\partial x^i\partial
  y^{q-i}}\right)\quad\text{otherwise,}\end{cases}
\]
so that  $(f,g)_q\in V_{n+m-2q}$. We may identify $\bigl(V,\langle
.\vert.\rangle\bigr)$ with $\bigl(V_1,(.,.)_1\bigr)$ (\cite[Lemma
2.2]{BDE}). Then (\cite[Theorem 3.2]{BDE}), the simple Lie algebra
of type $G_2$ appears as the $\bZ_2$-graded Lie algebra
\[
\frg_2=\Bigl(\frsp(V)\oplus V_2\bigr)\oplus \bigl(V\otimes
V_3\bigr),
\]
where $V_2$ is a Lie algebra with bracket $[f,g]=(f,g)_1$,
$\frsp(V)\oplus V_2$ is the even subalgebra of $\frg_2$, with
$[\frsp(V),V_2]=0$, $\frsp(V)$ acts naturally on $V\otimes V_3$ on
the first slot, $V_2$ acts on the second slot of $V\otimes V_3$ by
means of $[f,a\otimes F]=\tfrac{3}{2}a\otimes (f,F)_1$, for any
$f\in V_2$, $a\in V$ and $F\in V_3$, and finally,
\[
[a\otimes F,b\otimes G]=(F,G)_3\gamma_{a,b}+2\langle a\vert
b\rangle (F,G)_2\in \frsp(V)\oplus V_2,
\]
for any $a,b\in V$, $F,G\in V_3$. (Note that we have multiplied by
$-2$ the Lie bracket of odd elements in \cite[Theorem 3.2]{BDE}
and have used \cite[(2.5)]{BDE}.)

Theorem \ref{th:STSLie} now shows that $T=V_3$, with alternating
bilinear form given by $(.,.)_3$ and triple product
$[FGH]=3\bigl((F,G)_2,H)_1$, for any $F,G,H\in V_3$, is a \STS,
whose associated Lie algebra $\frg(T)$ is $\frg_2$ (simple of type
$G_2$). \qed

\end{examples}

\bigskip

The relationship of the \STSs\ described in \ref{ex:examples} with
the \STSs\ described in Theorem \ref{th:charne23} is given in the
next result:

\begin{theorem}\label{th:STSJordanLie}
Let $k$ be an algebraically closed field of characteristic $\ne
2,3$ and let $\T$ be a simple \STS. Then:
\begin{alphaenumerate}
\item If $[xyz]=\psi_{x,y}(z)$ for any $x,y,z\in T$ (item (i) in
\ref{th:charne23}), then $\frg(T)$ is isomorphic to the symplectic
Lie algebra $\frsp(V\perp T)$ as in Examples \ref{ex:examples}
(symplectic case).

\item If $\T$ is the \STS\ associated to a Jordan algebra $J$
which is either $0$ or the Jordan algebra of a nondegenerate
quadratic form with basepoint $J=\calJord(q,e)$, then $T$ is
isomorphic to the \STS\ in Examples \ref{ex:examples} (special
case), for a suitable vector space $W$ with $\dim W=1+\dim J$. In
particular, $\frg(T)$ is isomorphic to $\frpsl(V\oplus W)\simeq
\frpsl_{3+\dim J}(k)$.

\item If $\T$ is the \STS\ associated to the Jordan algebra of a
nondegenerate cubic form of type $J=k\times \calJord(q,e)$, then
$\T$ is isomorphic to the \STS\ in Examples \ref{ex:examples}
(orthogonal case), for a vector space $W$ endowed with a
nondegenerate symmetric bilinear form, with $\dim W=1+\dim J$. In
particular, $\frg(T)$ is isomorphic to the orthogonal Lie algebra
$\frso\bigl((V\otimes V)\oplus W\bigr)\simeq \frso_{5+\dim J}(k)$.

\item If $\T$ is the \STS\ associated to the Jordan algebra $J=k$
with cubic form $n(\alpha)=\alpha^3$, then $T$ is isomorphic to
the \STS\ defined on $V_3$ in Examples \ref{ex:examples}
(G$_2$-case). In particular, $\frg(T)$ is the simple Lie algebra
of type $G_2$.

\item If $\T$ is the \STS\ associated to the Jordan algebra
$H_3(C)$, where $C$ is either $k$, $k\times k$, $\Mat_2(k)$ or the
algebra of split octonions, then $\frg(T)$ is isomorphic,
respectively, to the simple Lie algebras of type $F_4$, $E_6$,
$E_7$ and $E_8$.
\end{alphaenumerate}
\end{theorem}
\begin{proof}
Part (a) is clear from the arguments in Example \ref{ex:examples}
(symplectic case) .

Also, if either $J=0$ or $J=\calJord(q,e)$ for a nondegenerate
quadratic form $q$ and basepoint $e$, then one has $T=W\oplus
\tilde W$ with $W=\left\{\bigl(\begin{smallmatrix} \alpha &a\\
0&0\end{smallmatrix}\bigr): \alpha\in k,\ a\in
J\right\}$ and $\tilde W= \left\{\bigl(\begin{smallmatrix} 0&0\\
b&\beta\end{smallmatrix}\bigr): \beta\in k,\ b\in J\right\}$. The
subspaces $W$ and $\tilde W$ are maximal isotropic subspaces of
$T$ and we will identify $\tilde W$ with the dual $W^*$ by means
of $x_2\in \tilde W\mapsto (x_2\vert .):W\rightarrow k$. Moreover,
the expression of the triple product in \eqref{eq:FJJFbis} gives
$[WWT]=0=[W^*W^*T]$. Besides, for any $x=\bigl(\begin{smallmatrix}
\alpha &a\\ 0&0\end{smallmatrix}\bigr)\in W$, $f=\bigl(\begin{smallmatrix} 0&0\\
b&\beta\end{smallmatrix}\bigr)\in W^*$ and
$y=\bigl(\begin{smallmatrix} \gamma &c\\
0&0\end{smallmatrix}\bigr)\in W$,
\[
\begin{split}
[xfy]&=\begin{pmatrix}
\bigl(-3\alpha\beta+t(a,b)\bigr)\gamma+2\alpha t(b,c)&
  \bigl(-\alpha\beta +t(a,b)\bigr)c+2\bigl(t(b,c)-\beta\gamma\bigr)a\\
  0&0\end{pmatrix}\\
  &=-(x\vert f)y-2(y\vert f)x,
\end{split}
\]
that is, equation \eqref{eq:xfy} is satisfied for any $x,y\in W$
and $f\in W^*$. Also, since $([xfy]\vert g)+(y\vert [xfg])=0$ for
any $x,y\in W$ and $f,g\in W^*$, or by direct computation as
above, it follows that equation \eqref{eq:xfg} is satisfied too.
Hence part (b) follows.

Now, in item (ii) of Theorem \ref{th:charne23}, the dimension of
$T$ is always even. If $\dim T=4$, then either $T$ is associated
to the unique Jordan algebra $J=\calJord(q,e)$ of a nondegenerate
quadratic form with basepoint of dimension $1$, which is already
known to be isomorphic to the \STS\ in Examples \ref{ex:examples}
(orthogonal case) with $\dim W=2$; or $J=k$ with cubic form
$n(\alpha)=\alpha^3$ and basepoint $1$. Then, necessarily, this
case corresponds to the G$_2$-case in Examples \ref{ex:examples}.

On the other hand, if $\dim T>4$, then either $T$ is associated to
the Jordan algebra of hermitian matrices $H_3(C)$, and hence the
associated Lie algebra $\frg(T)$ is $F_4,E_6,E_7$ or $E_8$ (this
goes back to Freudenthal \cite{FrII,FrVIII}, see also
\cite{YamAs}); or to the unique Jordan algebra $J=k\oplus
\calJord(q,e)$ with $\dim J=\tfrac{1}{2}(\dim T-1)$. The only
possibility left here is case (c).
\end{proof}

\begin{remark}\label{re:STSexceptional}
For a different construction of the exceptional \STSs\ in item (e)
of Theorem \ref{th:STSJordanLie}, see \cite{MSII}.
\end{remark}

\bigskip

Note that the \STSs\ that appear in Theorem \ref{th:charne23} make
sense too over fields of characteristic $3$, with the exceptions
in case (ii) of $J=0$, because then $[...]$ is the trivial product
by equations \eqref{eq:xfy} and \eqref{eq:xfg}, or $J=k$ with
$n(\alpha)=\alpha^3$ for any $\alpha\in k$, because then
$t(\alpha)=3\alpha=0$ for any $\alpha\in k$, so that the trace
form is trivial (this has to do with the fact that the exceptional
Lie algebra of type $G_2$ is no longer simple in characteristic
$3$ --see, for instance, \cite{Albercaetal}--). Theorem
\ref{th:STSFTS} and Brown's classification of the Freudenthal
triple systems in characteristic $3$ \cite{BrownFTS} yield:

\begin{theorem}\label{th:char3}
Let $k$ be an algebraically closed field of characteristic $3$ and
let $\T$ be a simple \STS, then either:
\begin{romanenumerate}
\item $[xyz]=\psi_{x,y}(z)$ for any $x,y,z\in T$; or
\item there exists a Jordan algebra $J$ of either a nondegenerate
quadratic form with basepoint, or of a nondegenerate cubic form
 with basepoint, such that $T$ is given
by the formulas in \eqref{eq:FJJF} and \eqref{eq:FJJFbis}; or
\item $\dim T=2$, and hence $T$ is described in item \textrm{(ii)} of
Proposition \ref{pr:STS2}; or
\item up to isomorphism, $T$ is the \STS\ described in Proposition
\ref{pr:STS8}.
\end{romanenumerate}
\end{theorem}
\begin{proof}
According to \cite{BrownFTS}, the only differences in
characteristic $3$ are given by the two dimensional \STSs\ and a
unique exceptional \STS\ of dimension $8$ which, therefore, must
be the one in Proposition \ref{pr:STS8}. Note that if $T$ is this
\STS, $\inder(T)$ is the Lie algebra $L(1)$ of Kostrikin, so this
\STS\ is indeed different from the previous ones.
\end{proof}

\bigskip

\begin{remark}\label{re:STSLiechar3}
The associated Lie algebras $\frg(T)$ of the algebras in item
(iii) of Theorem \ref{th:char3} are given in equation
\eqref{eq:gt2e}, Proposition \ref{pr:Lepsilon} and Remark
\ref{re:Lepsilon}. On the other hand, the simple Lie algebra
$\frg(T)$ associated to the eight-dimensional \STS\ in item (iv)
of Theorem \ref{th:char3} is a simple Lie algebra of dimension
$29$, which is specific of characteristic $3$ (see
\cite{BrownCont,Brown29}). This Lie algebra is the contragredient
Lie algebra with Cartan matrix $\left(\begin{smallmatrix}2&-1&0\\ -1&2&-1\\
0&-1&0\end{smallmatrix}\right)$.

As for the Lie algebras associated to the \STSs\ in items (i) and
(ii) of Theorem \ref{th:char3}, everything works as for
characteristic $\ne 2,3$, with the only exception of the \STS\
associated to the Jordan algebra $J=H_3(k\times k)=\Mat_3(K)^+$.
What happens here is that the split Lie algebra of type $E_6$ in
characteristic $3$ (obtained by taking a Chevalley basis of the
corresponding complex Lie algebra, then the $\bZ$-algebra spanned
by this basis, and finally tensoring with the field) is no longer
simple, but has a one-dimensional center (see, for instance,
\cite[\S 3]{VK}). Modulo this center, one gets a simple Lie
algebra of dimension $77$, which is precisely the Lie algebra
$\frg(T)$ in this case. This can be deduced from \cite{MSII}, we
will not go into details.
\end{remark}

\section{Lie superalgebras and symplectic triple
systems}\label{se:STSLie}

An interesting feature in characteristic $3$ is that any
symplectic triple system is an anti-Lie triple system, that is,
the odd part of a Lie superalgebra. More precisely, we can forget
about $\frsp(V)$ and $V$ in Theorem \ref{th:STSLie} and get, not a
$\bZ_2$-graded Lie algebra, but a Lie superalgebra. The precise
statement is the following:

\begin{theorem}\label{th:STSsuperLie}
Let $\T$ be a symplectic triple system over a field of
characteristic $3$. Define the superalgebra
$\frtg=\frtg(T)=\frtg\subo\oplus\frtg\subuno$, with:
\[
\frtg\subo=\inder(T),\qquad\frtg\subuno=T,
\]
and superanticommutative multiplication given by:
\begin{itemize}
\item
$\frtg\subo$ is a Lie subalgebra of $\frtg$;
\item
$\frtg\subo$ acts naturally on $\frtg\subuno$, that is,
$[d,x]=d(x)$ for any $d\in\inder(T)$ and $x\in T$;
\item
$[x,y]=d_{x,y}=[xy.]$, for any $x,y\in T$.
\end{itemize}
Then $\frtg(T)$ is a Lie superalgebra. Moreover, $T$ is simple if
and only if so is $\frtg(T)$.
\end{theorem}
\begin{proof}
For any $x,y,z\in T$, using the symmetry in the first two
variables of $[...]$ \eqref{eq:STSa}, the skew-symmetry of
$(.\vert.)$ and \eqref{eq:STSb}, one gets:
\[
\begin{split}
[[x,y],z]+&[[y,z],x]+[[z,x],y]\\
 &=[xyz]+[yzx]+[zxy]\\
 &=[xyz]+[yzx]-2[zxy]\\
 &=\bigl([xyz]-[xzy]\bigr)+\bigl([zyx]-[zxy]\bigr)\\
 &=\Bigl((x\vert z)y-(x\vert y)z+2(y\vert z)x\Bigr)
    +\Bigl((z\vert x)y-(z\vert y)x+2(y\vert x)z\Bigr)\\
 &=3\Bigl((y\vert z)x-(x\vert y)z\Bigr)=0.
\end{split}
\]
This shows that $\frtg(T)$ is a Lie superalgebra.

Moreover, any ideal $I$ of $T$ induces the ideal $d_{I,T}\oplus I$
of $\frtg(T)$. Conversely, if $\fra=\fra\subo\oplus\fra\subuno$ is
a nonzero ideal of $\frtg(T)$, then $I=\fra\subuno\ne 0$, as
$\frtg\subo$ acts faithfully on $\frtg\subuno$. Then
$[TTI]=[\frtg\subo,\fra\subuno]\subseteq \fra\subuno=I$, while
$[ITT]=[[\fra\subuno,\frtg\subuno],\frtg\subuno]\subseteq
\fra\subuno=I$. Hence $I$ is a nonzero ideal of $T$. This proves
the last part of the Theorem.
\end{proof}

\medskip

Let $k$ be a field of characteristic $3$ and consider the \STSs\
in Examples \ref{ex:examples}.

\smallskip

In the \emph{symplectic case}, $[xyz]=\psi_{x,y}(z)$ for any
$x,y,z\in T$, so that $\inder(T)=\frsp(T)$ and
$\frtg(T)=\frsp(T)\oplus T$, which is isomorphic to the
orthosymplectic Lie superalgebra $\frosp(W)$ of the superspace
$W=W\subo\oplus W\subuno$, with $W\subo=ke$ (a vector space of
dimension $1$), $W\subuno=T$ and supersymmetric bilinear form
given by extending $(.\vert.)$ on $T$ by means of $(e\vert e)=1$
and $(e\vert T)=0$.

\smallskip

In the \emph{special case}: $T=W\oplus W^*$ with triple product
determined in \eqref{eq:xfyg}. The restriction map
$\inder(T)\rightarrow \frgl(W)$, $d\mapsto d\vert_W$, gives an
isomorphism between $\inder(T)$ and either $\frsl(W)$ or
$\frgl(W)$ (depending on the dimension of $W$ being or not
congruent to $1$ modulo $3$). Let $\tilde W$ be the superspace
with $\tilde W\subo=ke$ and $\tilde W\subuno=W$, and identify the
special linear superalgebra $\frsl(\tilde W)$ with the vector
space of $2\times 2$ matrices $\bigl\{\bigl(\begin{smallmatrix}
\trace A&f\\ x&A\end{smallmatrix}\bigr) :x\in W,\, f\in W^*,\,
A\in \frgl(W)\bigr\}$. Let $\tilde\frz$ be the center of
$\frsl(\tilde W)$ ($\tilde\frz=0$ if $\dim W\not\equiv 1$ modulo
$3$, and $\dim\tilde\frz=1$ otherwise), and let $\frpsl(\tilde
W)=\frsl(\tilde W)/\tilde\frz$ be the corresponding projective
special linear Lie superalgebra. Then the linear map
$\tilde\Phi:\frtg(T)\rightarrow \frpsl(\tilde W)$, such that:
\[
\begin{split}
\tilde\Phi\bigl((x,f)\bigr)&=\begin{pmatrix} 0&-f\\ x&0\end{pmatrix} +\tilde\frz\,,\\
\tilde\Phi(d)&=\begin{cases} \begin{pmatrix} 0&0\\
0&d\vert_W\end{pmatrix}+\tilde\frz\qquad\text{if $\dim W\equiv 1\pmod 3$,}\\
\noalign{\smallskip}
\begin{pmatrix} \tfrac{\trace d\vert_W}{1-\dim W}&0\\
 0&d\vert_W+\tfrac{\trace d\vert_W}{1-\dim W}id\end{pmatrix}
 +\tilde\frz\qquad\text{otherwise},
\end{cases}
\end{split}
\]
for any $x\in W$, $f\in W^*$ and $d\in \inder(T)$, is an
isomorphism of Lie superalgebras.

\smallskip

In the \emph{orthogonal case}, $T=V\otimes W$, where $V$ is a two
dimensional vector space endowed with a nonzero alternating
bilinear form $\langle.\vert.\rangle$, and $W$ is a vector space
of dimension $\geq 3$ endowed with a nondegenerate symmetric
bilinear form $q_W$. In this case, from Examples \ref{ex:examples}
(orthogonal case):
\[
\frtg(T)=\Bigl(\frsp(V)\oplus \frso(W,q_W)\Bigr)\oplus
\bigl(V\otimes W\bigr),
\]
with \eqref{eq:STSorthogonal} yielding
\[
[u\otimes x,v\otimes y]=\frac{1}{2}q_w(x,y)\gamma_{u,v}+\langle
u\vert v\rangle \sigma_{x,y},
\]
for any $u,v\in V$ and $x,y\in W$. Then $\frtg(T)$ is easily
checked to be isomorphic to the orthosymplectic Lie superalgebra
of the vector superspace $\tilde W$, with $\tilde W\subo=W$ and
$\tilde W\subuno =V$, endowed with the nondegenerate
supersymmetric bilinear form given by $q_W$ on $\tilde W\subo$ and
$\langle .\vert .\rangle $ on $\tilde W\subuno$.

\smallskip

Also, if $T$ is a simple \STS\ of dimension $2$, then by Theorem
\ref{pr:STS2} there is a scalar $0\ne \epsilon\in k$ such that
$T\simeq T_{2,\epsilon}$. From equation \eqref{eq:gt2e}, we
conclude that $\frtg(T)\subo\simeq \frsl_2(k)$ and
$\frtg(T)\subuno$ is the two dimensional irreducible module for
$\frtg(T)\subo$. Hence $\frtg(T)$ is isomorphic to the
orthosymplectic Lie superalgebra $\frosp_{1,2}(k)$.

\medskip

The conclusion is that only well-known simple Lie superalgebras
appear as $\frtg(T)$ for these `classical' \STSs. However, the Lie
superalgebras $\frtg(T)$ for the remaining simple \STSs\ in
Theorem \ref{th:char3} have no counterpart in characteristic $0$
(see \cite{Kac,Scheunert}), thus providing new simple Lie
superalgebras in characteristic $3$:

\begin{theorem}\label{th:newchar3}
Let $k$ be an algebraically closed field of characteristic $3$.
Then there are simple finite dimensional Lie superalgebras $\frtg$
over $k$ satisfying:
\begin{romanenumerate}
\item $\dim\frtg=18\,(=10+8)$, $\frtg\subo$ is the Kaplansky Lie
superalgebra $L(1)$ and $\frtg\subuno$ is its $8$-dimensional
irreducible module in Proposition \ref{pr:STS8}.
\item
$\dim\frtg=35\,(=21+14)$, $\frtg\subo$ is the symplectic Lie
algebra $\frsp_6(k)$ and $\frtg\subuno$ is a $14$-dimensional
irreducible module for $\frtg\subo$.
\item
$\dim\frtg=54\,(=34+20)$, $\frtg\subo$ is the projective special
Lie algebra $\frpsl_6(k)$ and $\frtg\subuno$ is a $20$-dimensional
irreducible module for $\frtg\subo$.
\item
$\dim\frtg=98\,(=66+32)$, $\frtg\subo$ is the orthogonal Lie
algebra $\frso_{12}(k)$ and $\frtg\subuno$ is a $32$-dimensional
irreducible module for $\frtg\subo$ (spin module).
\item
$\dim\frtg=189\,(=133+56)$, $\frtg\subo$ is the simple Lie algebra
of type $E_7$ and $\frtg\subuno$ is a $56$-dimensional irreducible
module for $\frtg\subo$.
\end{romanenumerate}
\end{theorem}
\begin{proof}
It is enough to consider the simple \STSs\ in Theorem
\ref{th:char3} not considered in the previous remarks. These are
the $8$-dimensional simple \STS\ in Proposition \ref{pr:STS8}, and
the simple \STSs\ corresponding to the Jordan algebras $J=H_3(C)$
for $C=k$, $k\times k$, $\Mat_2(k)$ or the split octonions. The
irreducibility in (i) follows from the description before
Proposition \ref{pr:STS8}, and for the remaining four cases it can
be checked from the description of these \STSs\ in \cite{MSII}.
\end{proof}

\medskip

\begin{remark}\label{re:newchar3}
For the Lie superalgebras in items (ii)--(v), the odd part
$\frtg\subuno$ is given by a simple \STS\ $T$, and the action of
$\frtg\subo=\inder(T)$ on $\frtg\subuno =T$ coincides with the
action inside the Lie algebra $\frg(T)$ in Theorem
\ref{th:STSLie}. Since in all these cases, $\frg(T)$ is a
restricted Lie algebra, it follows that $\frtg\subuno$ is a
restricted irreducible $\frtg\subo$-module. Actually, with some
care (using for instance \cite[\S 4]{MSII}), it is checked that,
with the ordering of the simple roots as in \cite{Bou}
$\frtg\subuno$ is the irreducible restricted $\frtg\subo$-module
of highest weight $\omega_3$ in cases (ii) and (iii), any of
$\omega_5$ or $\omega_6$ in case (iv) and $\omega_7$ in case (v).
(See also \cite{Krutelevich}.)
\end{remark}

\medskip

It must be remarked here that, in the proof of Theorem
\ref{th:STSsuperLie}, the fact that the alternating bilinear form
$(.\vert .)$ on $T$ is nonzero has played no role. This suggests
the next definition.

\begin{definition}\label{df:NullSTS}
Let $T$ be a vector space over a field $k$ endowed with a triple
product $T\times T\times T\rightarrow T$, $(x,y,z)\mapsto [xyz]$.
Then $(T,[...])$ is said to be a \emph{null symplectic triple
system} if it satisfies the following identities:
\begin{subequations}\label{eq:NullSTS}
\begin{align}
&[xyz]=[yxz]=[xzy],\label{eq:NullSTSa} \\
&[xy[uvw]]=[[xyu]vw]+[u[xyv]w]+[uv[xyw]].\label{eq:NullSTSb}
\end{align}
\end{subequations}
\end{definition}

\smallskip

That is, the triple product is symmetric on its arguments, and for
each $x,y\in T$, the linear $d_{x,y}=[xy.]$ is a derivation.
Again, the linear span of these maps $d_{x,y}$'s will be denoted
by $\inder(T)$.

For Freudenthal triple systems, Kantor also considered the
possibility of the alternating form to be zero in \cite{KantorFr}.

Note that the construction in Theorem \ref{th:STSLie} may be
modified for null \STSs\ to give, with the same kind of arguments:

\begin{proposition}\label{pr:NullSTSLie}
Let $\bigl(T,[...]\bigr)$ be a null symplectic triple system and
let $\bigl(V,\langle.\vert.\rangle\bigr)$ be a two dimensional
vector space endowed with a nonzero alternating bilinear form.
Define the $\bZ_2$-graded algebra
$\frg=\frg(T)=\frgo\oplus\frguno$ with
\[
\begin{cases}
\frgo= \inder(T)\,,\\
\frguno=V\otimes T\,,
\end{cases}
\]
and anticommutative multiplication given by:
\begin{itemize}
\item
$\frgo$ is a Lie subalgebra of $\frg$;
\item
$\frgo$ acts naturally on $\frguno$, that is, $[d,v\otimes
x]=v\otimes d(x)$ for any  $d\in\inder(T)$, $v\in V$, and $x\in
T$;
\item
for any $u,v\in V$ and $x,y\in T$:
\begin{equation}\label{eq:nulloddbracket}
[u\otimes x,v\otimes y]=
  \langle u\vert v\rangle d_{x,y}
\end{equation}
where $d_{x,y}=[xy.]$.
\end{itemize}
Then $\frg(T)$ is a Lie algebra. Moreover, $\frg(T)$ is simple if
and only if so is $\T$.

Conversely, given a $\bZ_2$-graded Lie algebra
$\frg=\frgo\oplus\frguno$, with $\frgo=\frs$ and $\frguno=V\otimes
T$, where $T$ is a module for $\frs$ and the multiplication of odd
elements is given by \eqref{eq:nulloddbracket} for a
skew-symmetric bilinear map $d_{.,.}:T\times T\rightarrow \frs$,
$(x,y)\mapsto d_{x,y}$; then $T$ is a null \STS\ with the triple
product defined by $[xyz]=d_{x,y}(z)$, for any $x,y,z\in T$.
\end{proposition}
\begin{proof}
The fact that $\frg(T)$ is a Lie algebra and the converse follow
easily. Now, if $I$ is an ideal of $T$ (note that for null \STSs\
the ideals are precisely the $\inder(T)$-submodules of $T$), then
$d_{I,T}\oplus (V\otimes I)$ is an ideal of $\frg(T)$. Conversely,
assume that $T$ is simple, so $\frguno$ is the direct sum of two
isomorphic irreducible modules for $\frgo$ and hence it cannot be
the adjoint module for $\frgo$. Thus, it is enough to prove that
$\frg$ has no proper homogeneous ideals. As in Theorem
\ref{th:STSLie}, we may assume that if
$\fra=\fra\subo\oplus\fra\subuno$ is a proper ideal of $\frg$,
then $0\ne \fra\subuno$; so by irreducibility of $T$,
$\fra\subuno=u\otimes T$ for some $0\ne u\in V$. Since
$\fra\subo\supseteq [\frguno,u\otimes T]=\inder(T)=\frgo$, and
then $\frguno=[\frgo,\frguno]\subseteq \fra$, a contradiction.
\end{proof}

\medskip

\begin{remark}\label{re:NullSTSLie}
Given a null \STS\ $(T,[...])$, consider the Lie algebra $\frg(T)$
defined in Proposition \ref{pr:NullSTSLie}. If $\{v,w\}$ is a
symplectic basis of $\bigl(V,\langle .\vert .\rangle\bigr)$, then
the decomposition $\frg(T)=(v\otimes T)\oplus
\inder(T)\oplus(w\otimes T)$ is a $3$-grading of $\frg(T)$. This
shows that the Lie algebras $\frg(T)$, for null \STSs, are
precisely the $3$-graded Lie algebras
$\frg=\frg_{-1}\oplus\frg_0\oplus\frg_1$ such that $\frg_1$ and
$\frg_{-1}$ are isomorphic as modules over $\frg_0$.
\end{remark}

\medskip

Over fields of characteristic $\ne 3$, there are no simple null
\STSs:

\begin{theorem}\label{th:NoNullSTS}
Let $k$ be a field of characteristic $\ne 2,3$. Then there are no
simple null \STSs\ over $k$.
\end{theorem}
\begin{proof}
Assume, on the contrary, that $\bigl(T,[...]\bigr)$ is a simple
null \STS\ over $k$. We will get a contradiction following several
steps:

\smallskip

\noindent\textbf{(i)} Let us prove first that for any $x\in T$,
$d_{x,x}^2=0$:

In fact, from \eqref{eq:NullSTS} we conclude that for any $x,y\in
T$, $[d_{x,y},d_{x,x}]=2d_{[xyx],x}$. But also,
$[d_{x,y},d_{x,x}]=-[d_{x,x},d_{x,y}]=-d_{[xxx],y}-d_{x,[xxy]}$.
Therefore
\begin{equation}\label{eq:3dxxxy}
3d_{x,[xxy]}=-d_{[xxx],y},
\end{equation}
for any $x,y\in T$. If this is applied to $x$, one gets by
\eqref{eq:NullSTSa} that $3[xx[xxy]]=-[xy[xxx]]=-3[xx[xxy]]$
(since $d_{x,y}$ is a derivation). Therefore, $6[xx[xxy]]=0$ and
hence $d_{x,x}^2=0$.

By linearization of $d_{x,x}^2=0$ we get
\begin{equation}\label{eq:dxxdyy}
d_{x,x}d_{y,y}+d_{y,y}d_{x,x}+4d_{x,y}^2=0
\end{equation}
for any $x,y\in T$.

\smallskip

\noindent\textbf{(ii)} For any $x,y\in T$,
$d_{x,x}d_{y,y}d_{x,x}=-d_{[xxy],[xxy]}$:

Since $d_{x,x}^2=0$,
$[d_{x,x},[d_{x,x},d_{y,y}]]=-2d_{x,x}d_{y,y}d_{x,x}$. But, by the
derivation property and the symmetry of the triple product, also
$[d_{x,x},[d_{x,x},d_{y,y}]]=2[d_{x,x},d_{[xxy],y}]=2d_{[xxy],[xxy]}$,
whence the result. As a consequence, we get that for any $x\in T$:
\begin{equation}\label{eq:xxTxxTT}
[[xxT][xxT]T]\subseteq [xxT].
\end{equation}

\smallskip

\noindent\textbf{(iii)} There are elements $x\in T$ such that
$d_{x,x}=0$:

Take $0\ne S$ a minimal nonzero subspace of $T$ with the property
that $[SST]\subseteq S$. By equation \eqref{eq:xxTxxTT}, for any
$0\ne x\in S$, $S_x=[xxT]\subseteq S$ and also $[S_xS_xT]\subseteq
S_x$. Therefore, either $d_{x,x}=0$ and we are done, or $S_x=S$
and hence there is an element $y\in S$ such that $[xxy]=x$. But
then $[xxx]=d_{x,x}(x)=d_{x,x}^2(y)=0$, so by \eqref{eq:3dxxxy}
$3d_{x,x}=3d_{x,[xxy]}=0$, and thus $x=d_{x,x}(y)=0$, a
contradiction.

\smallskip

\noindent\textbf{(iv)} $T$ is the linear span of the elements $x$
with $d_{x,x}=0$:

By \textbf{(ii)}, if $d_{x,x}=0$, then also $d_{[yyx],[yyx]}=0$.
Thus the span of $\{x\in T : d_{x,x}=0\}$ is
$\inder(T)$-invariant, and hence an ideal of $T$.

\smallskip

\noindent\textbf{(v)} A subtriple $S$ of $T$ is said to be
\emph{solvable} if $0$ belongs to the chain of subtriples defined
by $S^{(0)}=S$ and $S^{(i+1)}=[S^{(i)}S^{(i)}S^{(i)}]$ for any
$i\geq 0$. Then, if $S$ is a solvable subtriple of $T$, the
subalgebra of the associative algebra $\End_k(T)$ generated by
$d_{S,S}$ is nilpotent:

This can be proved by induction on the minimal natural number $n$
such that $S^{(n)}=0$. This is trivial for $n=0$. Otherwise, let
$R$ be the subalgebra of $\End_k(T)$ generated by $d_{S,S}$ and
let $A$ be the subalgebra generated by $d_{S^{(1)},S^{(1)}}$. By
induction hypothesis, $A$ is nilpotent. Consider also the
subalgebra $B$ of $\End_k(T)$ generated by $d_{S,S^{(1)}}$. Let us
first check that $B$ is nilpotent.

Since $[d_{S,S},d_{S,S^{(1)}}]\subseteq d_{S^{(1)},S^{(1)}}$, it
follows that $RA\subseteq AR+A$, and hence  $I=AR+A$ is a
nilpotent ideal of $R$. Now, for any $s,t\in S$ and $s_1,t_1\in
S^{(1)}$, $[d_{s,s_1},d_{t,t_1}]=d_{[ss_1t],t_1}+d_{t,[ss_1t_1]}$
and $d_{[ss_1t],t_1}$ belongs to $d_{S^{(1)},S^{(1)}}\subseteq
A\subseteq B\cap I$. Therefore, the generating subspace
$\bigl(d_{S,S^{(1)}}+ B\cap I\bigr)/B\cap I$ of $B/B\cap I$ is
closed under commutation. Also, for any $x\in S$ and $y\in
S^{(1)}$,
$d_{x,y}^2=-\tfrac{1}{4}\bigl(d_{x,x}d_{y,y}+d_{y,y}d_{x,x}\bigr)\in
RA+AR\subseteq I$ (by \eqref{eq:dxxdyy}). Hence, by Engel's
Theorem \cite[Theorem II.1]{JacobsonJordan}, the subalgebra
$B/B\cap I$ is nilpotent and, since $I$ is already known to be
nilpotent, we conclude that the subalgebra $B$ is nilpotent.

On the other hand, $[d_{S,S},d_{S,S^{(1)}}]\subseteq
d_{S,S^{(1)}}$, so that $RB\subseteq BR+B$. Hence $J=BR+B$ is a
nilpotent ideal of $R$. Also, $[d_{S,S},d_{S,S}]\subseteq
d_{S,S^{(1)}}\subseteq B\subseteq J$. Thus, $R/J$ is a commutative
algebra generated by nilpotent elements, and hence $R/J$ is
nilpotent, and so is $R$, as required.

\smallskip

\noindent\textbf{(vi)} To get to a contradiction, let $S$ be a
maximal solvable subtriple of $T$. Since $T$ is simple, $S\ne T$.
By \textbf{(v)}, $d_{S,S}$ acts nilpotently on $T$, so that,
because of \textbf{(iv)}, there is an element $w\in T\setminus S$
such that $d_{w,w}=0$ and a natural number $n$ such that
$d_{S,S}^n(w)\not\in S$, but $d_{S,S}^{n+1}(w)\subseteq S$. Hence,
there are elements $x_1,\ldots,x_n\in S$ such that
$x=d_{x_1,x_1}\cdots d_{x_n,x_n}(w)\not\in S$,
$d_{S,S}(x)\subseteq S$ and, because $d_{w,w}=0$ and
\textbf{(ii)}, also $d_{x,x}=0$. Take $S'=S+kx$. Then
$[S'S'S']\subseteq [SSS]+[SSx]\subseteq S$. Hence $S'$ is solvable
and larger than $S$, a contradiction.
\end{proof}

\medskip

\begin{remark} This proof is inspired in similar proofs for Jordan
pairs \cite[Ch.\,14]{Loos}. Actually, observe that if $\{v,w\}$ is
any basis of $V$, then the decomposition $\frg(T)=(v\otimes T)
\oplus \inder T\oplus (w\otimes T)$ is a $3$-grading of $\frg(T)$,
and hence the pair $(v\otimes T,w\otimes T)$ is a Jordan pair
(characteristic $\ne 3$), which is simple if so is $T$. One can
use then known results on Jordan pairs to get Theorem
\ref{th:NoNullSTS}, but we have preferred to give a self contained
proof.
\end{remark}

\medskip

By Theorem \ref{th:STSFaulkner}, the null \STSs\ correspond
bijectively to the Faulkner ternary algebras with trivial
alternating bilinear form. Therefore, we get the following
improvement of \cite[Lemma 3.1]{FaulknerFerrar}, which in turn,
shows that \cite{FaulknerFerrar} contains a complete
classification of the simple Faulkner ternary algebras over
algebraically closed fields of characteristic $\ne 2,3$.

\begin{corollary}\label{co:FFcharne3}
Let $k$ be a field of characteristic $\ne 2,3$. Then a Faulkner
ternary algebra $\bigl(T,\langle...\rangle,(.\vert .)\bigr)$ over
$k$ is simple if and only if $(.\vert .)$ is nondegenerate.
\end{corollary}
\begin{proof}
By Theorem \ref{th:NoNullSTS}, if
$\bigl(T,\langle...\rangle,(.\vert .)\bigr)$ is simple, then
$(.\vert .)$ is not $0$ and now Proposition \ref{pr:STSsimple} (or
\cite[Lemma 3.1]{FaulknerFerrar}) applies.
\end{proof}

\medskip

The comments in the paragraph previous to Proposition
\ref{pr:STS2} show that any two dimensional simple \STS\ over a
field of characteristic $3$ is actually a simple null \STS. The
following conjecture sounds plausible.

\begin{conjecture} The dimension of any simple null \STS\ over
a field of characteristic $3$ is two.
\end{conjecture}

\section{Orthogonal triple systems}\label{se:OTS}

If the symmetry and skew-symmetry in the definition of a
symplectic triple system (Definition \ref{df:STS}) are
interchanged, one obtains the definition of orthogonal triple
systems, which first appeared in \cite[Section V]{OkuboI}:

\begin{definition}\label{df:OTS}
Let $T$ be a vector space endowed with a nonzero symmetric
bilinear form $(.\vert.):T\times T\rightarrow k$, and a triple
product $T\times T\times T\rightarrow T$: $(x,y,z)\mapsto [xyz]$.
Then $\bigl(T,[...],(.\vert.)\bigr)$ is said to be a
\emph{orthogonal triple system} if it satisfies the following
identitities:
\begin{subequations}\label{eq:OTS}
\begin{align}
&[xxy]=0\label{eq:OTSa}\\
&[xyy]=(x\vert y)y-(y\vert y)x\label{eq:OTSb}\\
&[xy[uvw]]=[[xyu]vw]+[u[xyv]w]+[uv[xyw]]\label{eq:OTSc}\\
&([xyu]\vert v)+(u\vert [xyv])=0\label{eq:OTSd}
\end{align}
\end{subequations}
for any elements $x,y,u,v,w\in T$.
\end{definition}

Note that \eqref{eq:OTSb} can be written as
\begin{equation}\label{eq:OTSbb}
[xyy]=\sigma_{x,y}(y)
\end{equation}
with $\sigma_{x,y}(z)=(x\vert z)y-(y\vert z)x$. (If $(.\vert.)$ is
nondegenerate, the maps $\sigma_{x,y}$ span the orthogonal Lie
algebra $\frso(T)$.)

Also, as for the symplectic case, \eqref{eq:OTSc} is equivalent to
$d_{x,y}=[xy.]$ being a derivation of the triple system. Let
$\inder(T)$ be the linear span of $\{d_{x,y}:x,y\in T\}$, which is
a Lie subalgebra of $\End(T)$. Equation \eqref{eq:OTSd} is
equivalent to $d_{x,y}\in \frso(T)$ for any $x,y\in T$.

Here condition \eqref{eq:OTSd} is always a consequence of the
previous ones. This is trivial if $\dim T=1$; otherwise take
linearly independent elements $u,v\in T$. Then \eqref{eq:OTSc}
with $w=v$ gives, using \eqref{eq:OTSb}, that
\[
\Bigl(([xyu]\vert v)+(u\vert [xyv])\Bigr)v=2([xyv]\vert v)u,
\]
whence the claim. Anyway, we have preferred to keep
\eqref{eq:OTSd} in the definition.

The definition of homomorphism between two \OTSs\ and of ideal and
simplicity are the natural ones.

Orthogonal triple systems are strongly related to another class of
triple systems: the $(-1,-1)$ balanced Freudenthal Kantor triple
systems (see \cite{EKO} and the references there in).

With the same kind of arguments, as in the proof of Proposition
\ref{pr:STSsimple} (see also \cite[Theorem 2.2]{EKO}) one gets:

\begin{proposition}\label{pr:OTSsimple}
Let $\T$ be an \OTS. Then $\T$ is simple if and only if the
bilinear form $(.\vert.)$ is nondegenerate.
\end{proposition}

And by mimicking the proof of Theorem \ref{th:STSLie}, or using
\cite[Theorem 2.1]{EKO}, it is shown how \OTSs\ are related to a
specific class of Lie superalgebras:

\begin{theorem}\label{th:OTSsuperLie}
Let $\bigl(T,[...],(.\vert .)\bigr)$ be an \OTS\  and let
$\bigl(V,\langle.\vert.\rangle\bigr)$ be a two dimensional vector
space endowed with a nonzero alternating bilinear form. Define the
superalgebra $\frg=\frg(T)=\frgo\oplus\frguno$ with
\[
\begin{cases}
\frgo= \frsp(V)\oplus \inder(T)&\text{(direct sum of ideals),}\\
\frguno=V\otimes T\,,
\end{cases}
\]
and superanticommutative multiplication given by:
\begin{itemize}
\item
$\frgo$ is a Lie subalgebra of $\frg$;
\item
$\frgo$ acts naturally on $\frguno$, that is,
\[
[s,v\otimes x]=s(v)\otimes x,\qquad [d,v\otimes x]=v\otimes d(x),
\]
for any $s\in \frsp(V)$, $d\in\inder(T)$, $v\in V$, and $x\in T$;
\item
for any $u,v\in V$ and $x,y\in T$:
\begin{equation}\label{eq:superoddbracket}
[u\otimes x,v\otimes y]=
  (x\vert y)\gamma_{u,v} -\langle u\vert v\rangle d_{x,y}
\end{equation}
where $\gamma_{u,v}=\langle u\vert .\rangle v+\langle v\vert
.\rangle u$ and $d_{x,y}=[xy.]$.
\end{itemize}
Then $\frg(T)$ is a Lie superalgebra. Moreover, $\frg(T)$ is
simple if and only if so is $\T$.

Conversely, given a  Lie superalgebra $\frg=\frgo\oplus\frguno$
with
\[
\begin{cases}
\frgo=\frsp(V)\oplus \frs &\text{(direct sum of ideals),}\\
\frguno=V\otimes T&\text{(as a module for $\frgo$),}
\end{cases}
\]
where $T$ is a module for $\frs$, by $\frgo$-invariance of the Lie
bracket, equation \eqref{eq:superoddbracket} is satisfied for a
symmetric bilinear form $(.\vert .):T\times T\rightarrow k$ and a
skew-symmetric bilinear map $d_{.,.}:T\times T\rightarrow \frs$.
Then, if $(.\vert .)$ is not $0$ and a triple product on $T$ is
defined by means of $[xyz]=d_{x,y}(z)$, $\bigl(T,[...],(.\vert
.)\bigr)$ is an \OTS.
\end{theorem}

\medskip

Over fields of characteristic $0$, the classification of the
simple $(-1,-1)$ balanced Freudenthal Kantor triple systems in
\cite[Theorem 4.3]{EKO} immediately implies the following
classification of the simple \OTSs.

\begin{theorem}\label{th:OTSchar0}
Let $k$ be a field of characteristic zero and let $\T$ be a simple
\OTS\ over $k$. Then either:
\begin{romanenumerate}
\item
 $[xyz]=(x\vert z)y-(y\vert z)x\,
\bigl(=\sigma_{x,y}(z)\bigr)$ for any $x,y,z\in T$
(\emph{orthogonal type}). In this case $\inder(T)=\frso(T)$.
\item
 There is a quadratic \'etale algebra $K$ over $k$ such that
$T$ is a free $K$-module of rank at least $3$, endowed with a
nondegenerate hermitian form $h:T\times T\rightarrow K$
($h(x,y)=\overline{h(y,x)}$, $h(rx,y)=rh(x,y)$, for any $x,y\in T$
and $r\in K$, where $r\mapsto \bar r$ is the standard involution
on $K$) such that
\begin{equation}\label{eq:OTSunitarian}
\begin{cases}
(x\vert x)=h(x,x)\\
[xyz]=h(z,x)y-h(z,y)x+\frac{1}{2}\bigl(h(x,y)-h(y,x)\bigr)z
\end{cases}
\end{equation}
for any $x,y,z\in T$ (\emph{unitarian type}). In this case
$\inder(T)$ is the unitarian Lie algebra $\fru(T,h)=\{ f\in
\End_K(T): h(f(x),y)+h(x,f(y))=0,\ \text{for any $x,y\in T$}\}$.
\item
 There is a quaternion algebra $Q$ over $k$ such that $T$ is a
 free left $Q$-module of rank $\geq 2$, endowed with a
 nondegenerate hermitian form $h:T\times T\rightarrow Q$
 satisfying \eqref{eq:OTSunitarian} (\emph{symplectic type}). In
 this case $\inder(T)$ is the direct sum of the symplectic Lie algebra
 $\frsp(T,h)=\{ f\in\End_Q(T): h(f(x),y)+h(x,f(y))=0,\
 \text{for any $x,y\in T$}\}$ and of the three dimensional simple Lie algebra
 $[Q,Q]$.
\item
 $\dim_kT=4$ and there is a nonzero skew-symmetric multilinear
 map $\Phi: T\times T\times T\times T\rightarrow k$ such that, for
 any $x,y,z,t\in T$,
 \begin{equation}\label{eq:OTS4}
 [xyz]=\{ xyz\}+\sigma_{x,y}(z),
 \end{equation}
 where $\{ ...\}$ is defined by means of $\Phi(x,y,z,t)=(\{
 xyz\}\vert t)$. In this case there is a scalar $0\ne \mu\in k$
 such that $(\{a_1a_2a_3\}\vert \{b_1b_2b_3\})=\mu\det\bigl((a_i\vert
 b_j)\bigr)$ (\emph{D$_\mu$-type}). For $\mu=1$, $\inder(T)$ is a
 three dimensional simple Lie algebra, while for $\mu\ne 0,1$,
 $\inder(T)=\frso(T)$.
\item
 $\dim_kT=7$ and there is an eight-dimensional Cayley algebra $C$
 over $k$ with trace $t$ and a nonzero scalar $\alpha\in k$ such
 that $T=C_0=\{x\in C: t(x)=0\}$ and for any $x,y,z\in T$,
\begin{equation}\label{eq:OTSG}
\begin{cases}
(x\vert y)=-2\alpha t(xy)\\
[xyz]=\alpha D_{x,y}(z)
\end{cases}
\end{equation}
where $D_{x,y}=[L_x,L_y]+[L_x,R_y]+[R_x,R_y]\in \der(C)$ ($L_x$
and $R_x$ denote the left and right multiplications by $x$ in $C$)
(\emph{G-type}). In this case $\inder(T)$ is the Lie algebra of
derivations  of $C$: $\der(C)$ (viewed as a Lie subalgebra of
$\frgl(C_0)$), which is a simple Lie algebra of type $G_2$.
\item
 $\dim_kT=8$ and $T$ is endowed with a $3$-fold vector cross
 product $X$ of type $I$ (see \cite{Elduque3fold} and the
 references there in) relative to a nondegenerate symmetric
 bilinear form $b(.,.)$ (so that
 \[
 b\bigl(X(a_1,a_2,a_3),X(a_1,a_2,a_3)\bigr)=\det\bigl(
 b(a_i,a_j)\bigr)
 \] for any $a_1,a_2,a_3\in T$)
 such that
 \begin{equation}\label{eq:OTSF}
 [xyz]=X(x,y,z)+3\tau_{x,y}(z)
 \end{equation}
 for any $x,y,z\in T$, where $\tau_{x,y}=b(x,.)y-b(y,.)x$ and
 $(.\vert .)=3b(.,.)$ (\emph{F-type}). In this case $\inder(T)$ is
 isomorphic to the orthogonal Lie algebra $\frso(W,b)$, where $W$
 is any regular seven dimensional subspace of $T$ relative to $b$.
 Besides, $T$ is the spin module for $\inder(T)$.
\end{romanenumerate}
\end{theorem}

\medskip

\begin{remark}\label{re:OTSchar0}
Two \OTSs\ in different items in Theorem \ref{th:OTSchar0} are
never isomorphic, and the conditions for isomorphism among \OTSs\
in the same item can be read off \cite[Theorem 4.3]{EKO}. Only the
F-type has been described in a slightly different way as it is in
\cite{EKO}, where our $X$ denotes what appears as $\tfrac{1}{3}X$
there.
\end{remark}

\medskip

As for \STSs, the definition of \OTSs\ makes sense even if the
symmetric bilinear form involved is trivial.

\begin{definition}\label{df:NullOTS}
Let $T$ be a vector space over a field $k$ endowed with a triple
product $T\times T\times T\rightarrow T$, $(x,y,z)\mapsto [xyz]$.
Then $(T,[...])$ is said to be a \emph{null orthogonal triple
system} if it satisfies the following identities:
\begin{subequations}\label{eq:NullOTS}
\begin{align}
&[xyy]=[yyx]=0,\label{eq:NullOTSa} \\
&[xy[uvw]]=[[xyu]vw]+[u[xyv]w]+[uv[xyw]].\label{eq:NullOTSb}
\end{align}
\end{subequations}
\end{definition}

\medskip

The following counterpart to Proposition \ref{pr:NullSTSLie}
follows with the same sort of arguments used there:

\begin{proposition}\label{pr:NullOTSLie}
Let $\bigl(T,[...]\bigr)$ be a null \OTS\ and let
$\bigl(V,\langle.\vert.\rangle\bigr)$ be a two dimensional vector
space endowed with a nonzero alternating bilinear form. Define the
superalgebra  $\frg=\frg(T)=\frgo\oplus\frguno$ with
\[
\begin{cases}
\frgo= \inder(T)\,,\\
\frguno=V\otimes T\,,
\end{cases}
\]
and superanticommutative multiplication given by:
\begin{itemize}
\item
$\frgo$ is a Lie subalgebra of $\frg$,
\item
$\frgo$ acts naturally on $\frguno$, that is, $[d,v\otimes
x]=v\otimes d(x)$ for any  $d\in\inder(T)$, $v\in V$, and $x\in
T$;
\item
for any $u,v\in V$ and $x,y\in T$:
\begin{equation}\label{eq:nulloddOTSbracket}
[u\otimes x,v\otimes y]=
  \langle u\vert v\rangle d_{x,y}
\end{equation}
where $d_{x,y}=[xy.]$.
\end{itemize}
Then $\frg(T)$ is a Lie superalgebra. Moreover, $\frg(T)$ is
simple if and only if so is $\T$.

Conversely, given a  Lie superalgebra $\frg=\frgo\oplus\frguno$,
with $\frgo=\frs$ and $\frguno=V\otimes T$, where $T$ is a module
for $\frs$ and the multiplication of odd elements is given by
\eqref{eq:nulloddOTSbracket} for a bilinear symmetric map
$d_{.,.}:T\times T\rightarrow \frs$, $(x,y)\mapsto d_{x,y}$; then
$T$ is a null \OTS\ with the triple product defined by
$[xyz]=d_{x,y}(z)$, for any $x,y,z\in T$.
\end{proposition}

\medskip

\begin{remark}\label{re:NullOTSLie}
As for \STSs, given a null \OTS\ $(T,[...])$, consider the Lie
superalgebra $\frg(T)$ defined in Proposition \ref{pr:NullOTSLie}.
If $\{v,w\}$ is a symplectic basis of $\bigl(V,\langle .\vert
.\rangle\bigr)$, then the decomposition $\frg(T)=(v\otimes
T)\oplus \inder(T)\oplus(w\otimes T)$ is a consistent $3$-grading
of $\frg(T)$ ($\frgo=\frg_0$ and $\frguno=\frg_{-1}\oplus\frg_1$).
This shows that the Lie superalgebras $\frg(T)$, for null \OTSs,
are precisely the consistently $3$-graded Lie superalgebras
$\frg=\frg_{-1}\oplus\frg_0\oplus\frg_1$ such that $\frg_1$ and
$\frg_{-1}$ are isomorphic as modules over $\frg_0$.
\end{remark}

\medskip

The classification in Theorem \ref{th:OTSchar0} can be completed
by means of the following:

\begin{theorem}\label{th:NullOTSchar0}
Let $\T$ be a simple null \OTS\ over a field $k$ of characteristic
$0$. Then $\dim_kT=4$ and there is a nondegenerate symmetric
bilinear form $(.\vert .)$ and a nonzero multilinear
skew-symmetric map $\Phi:T\times T\times T\times T\rightarrow k$
such that
\begin{equation}\label{eq:NullOTS4}
\Phi(x,y,z,t)=([xyz]\vert t)
\end{equation}
for any $x,y,z,t\in T$. Conversely, any triple system defined by
means of \eqref{eq:NullOTS4} is a simple null \OTS.

Moreover, given two such null \OTSs\ $(T_1,[...]_1)$ and
$(T_2,[...]_2)$ over $k$ with associated multilinear maps $\Phi_1$
and $\Phi_2$ and nondegenerate symmetric bilinear forms
$(.\vert.)_1$ and $(.\vert.)_2$, let $\mu_l$ ($l=1,2$) be the
nonzero scalars such that
\[
\bigl([x_1x_2x_3]\vert
[y_1y_2y_3]\bigr)_l=\mu_l\det\Bigl((x_i\vert y_j)_l\Bigr)
\]
for any $x_i,y_i\in T_l$. Then $(T_1,[...]_1)$ is isomorphic to
$(T_2,[...]_2)$ if and only if there is a similarity $\varphi:
\bigl(T_1,(.\vert.)_1\bigr)\rightarrow \bigl(T_2,(.\vert
.)_2\bigr)$ of norm $\alpha$ (that is,
$\bigl(\varphi(x)\vert\varphi(y)\bigr)_2=\alpha(x\vert y)_1$ for
any $x,y\in T_1$) such that $\mu_1=\mu_2\alpha^2$.
\end{theorem}
\begin{proof}
Assume first that $\dim_kT=4$ and that there are maps $\Phi$ and
$(.\vert.)$ such that \eqref{eq:NullOTS4} is satisfied. Then
$[...]$ is skew-symmetric \eqref{eq:NullOTSa}, since so is $\Phi$.
Then, if $\{e_1,e_2,e_3,e_4\}$ is any basis of $T$ with
$\Phi(e_1,e_2,e_3,e_4)=1$, and $\{e^1,e^2,e^3,e^4\}$ is its dual
basis relative to $(.\vert.)$, the map $d_{e_1,e_2}=[e_1e_2.]$
annihilates $e_1$ and $e_2$ and takes $e_3$ to $e^4$ and $e_4$ to
$-e^3$. This argument shows easily that the maps $d_{e_i,e_j}$
($1\leq i<j\leq 4$) are linearly independent, and hence $\dim
d_{T,T}=6$. On the other hand, because of \eqref{eq:NullOTS4},
$d_{T,T}\subseteq \frso(T)$ so, by dimension count,
$d_{T,T}=\frso(T)$, which is contained in the special linear Lie
algebra $\frsl(T)$. But $\Phi$ is invariant under $\frsl(T)$, so
both $\Phi$ and $(.\vert .)$ are invariant under $d_{T,T}$. This
shows that $[...]$ is invariant too under $d_{T,T}$, which is
equivalent to condition \eqref{eq:NullOTSb} and, therefore,
$(T,[...])$ is a null \OTS, which is irreducible as a $\inder(T)
=\frso(T)$-module, and hence simple.

Now, if $\varphi:(T_1,[...]_1)\rightarrow (T_2,[...]_2)$ is an
isomorphism between two of these null \OTSs, then for any
$a_1,a_2\in T_1$, $[\varphi(a_1)\varphi(a_2)\,.\,]_2=\varphi
[a_1a_2\,.\,]_1\varphi^{-1}$, so that $\varphi$ induces an
isomorphism $\tilde\varphi:\frso(T_1)\rightarrow \frso(T_2)$,
$d\mapsto \varphi d\varphi^{-1}$. But then $(.\vert.)_1$ and
$\bigl(\varphi(.)\vert\varphi(.)\bigr)_2$ are both invariant under
the action of $\frso(T_1)$. Therefore, there is a nonzero scalar
$\alpha\in k$ such that
$\bigl(\varphi(x)\vert\varphi(y)\bigr)_2=\alpha(x\vert y)_1$ for
any $x,y\in T_1$, and $\varphi$ is a similarity of norm $\alpha$.
Moreover, for any $x_1,x_2,x_3,y_1,y_2,y_3\in T_1$:
\[
\begin{split}
\alpha\mu_1\det\bigl((x_i\vert y_j)_1\bigr)
 &=\alpha\bigl([x_1x_2x_3]_1\vert [y_1y_2y_3]_1\bigr)_1\\
 &=\bigl(\varphi([x_1x_2x_3]_1)\vert\varphi([y_1y_2y_3]_1\bigr)_2\\
 &=\bigl([\varphi(x_1)\varphi(x_2)\varphi(x_3)]_2\vert
     [\varphi(y_1)\varphi(y_2)\varphi(y_3)]_2\bigr)_2\\
 &=\mu_2\det\bigl((\varphi(x_i)\vert\varphi(y_j))_2\bigr)\\
 &=\mu_2\alpha^3\det\bigl((x_i\vert y_j)_1\bigr)
\end{split}
\]
so that $\mu_1=\mu_2\alpha^2$. Conversely, assume that
$\psi:\bigl(T_1,(.\vert)_1\bigr)\rightarrow \bigl(T_2,(.\vert
.)_2\bigr)$ is a similarity of norm $\alpha$ such that
$\mu_1=\mu_2\alpha^2$. Then, since the skew-symmetric multilinear
map on $T_1$ given by
$\Phi_1'(x,y,z,t)=\Phi_2(\psi(x),\psi(y),\psi(z),\psi(t))$ is (as
$\dim_kT_1=4$) a scalar multiple of $\Phi_1$, there is a nonzero
scalar $\beta\in k$ such that
$\Phi_2(\psi(x),\psi(y),\psi(z),\psi(t))=\beta\Phi_1(x,y,z,t)$ for
any $x,y,z,t\in T_1$, and hence we conclude that
\[
[\psi(x)\psi(y)\psi(z)]_2=\frac{\beta}{\alpha}\psi([xyz]_1)
\]
for any $x,y,z\in T_1$. But now, for any
$x_1,x_2,x_3,y_1,y_2,y_3\in T_1$:
\[
\begin{split}
\alpha\mu_1\det\bigl((x_i\vert y_j)_1\bigr)
 &=\alpha\bigl([x_1x_2x_3]_1\vert [y_1y_2y_3]_1\bigr)_1\\
 &=\bigl(\psi([x_1x_2x_3]_1)\vert\psi([y_1y_2y_3]_1)\bigr)_2\\
 &=\left(\tfrac{\alpha}{\beta}\right)^2\bigl([\psi(x_1)\psi(x_2)\psi(x_3)]_2\vert
     [\psi(y_1)\psi(y_2)\psi(y_3)]_2\bigr)_2\\
 &=\left(\tfrac{\alpha}{\beta}\right)^2\mu_2\det\bigl((\psi(x_i)\vert\psi(y_j))_2\bigr)\\
 &=\left(\tfrac{\alpha}{\beta}\right)^2\mu_2\alpha^3\det\bigl((x_i\vert y_j)_1\bigr)
\end{split}
\]
and, since $\mu_1=\mu_2\alpha^2$, we obtain
$\left(\tfrac{\alpha}{\beta}\right)^2=1$ or $\beta=\pm\alpha$. In
case $\beta=\alpha$, $\psi$ is an isomorphism between
$(T_1,[...]_1)$ and $(T_2,[...]_2)$, while if $\beta=-\alpha$,
$\psi$ satisfies $\psi([xyz]_1)=-[\psi(x)\psi(y)\psi(z)]_2$ for
any $x,y,z\in T_1$. But we can always take an orthogonal map
$\sigma$ of $\bigl(T_1,(.\vert.)_1\bigr)$ (that is
$(\sigma(x)\vert \sigma(y))_1=(x\vert y)_1$ for any $x,y$) with
$\det\sigma=-1$. Then, for any $x,y,z,t\in T_1$,
$\Phi_1(\sigma(x),\sigma(y),\sigma(z),\sigma(t))=-\Phi_1(x,y,z,t)$,
so that $\sigma([xyz]_1)=-[\sigma(x)\sigma(y)\sigma(z)]_1$ and the
composite map $\psi\sigma$ is an isomorphism.

Finally, let us show that these four dimensional simple null
\OTSs\ are the only ones. To do so, we may assume that $k$ is
algebraically closed. Note that, by Proposition
\ref{pr:NullOTSLie}, any simple null \OTS\ $(T,[...])$ gives rise
to a simple Lie superalgebra $\frg=\frg(T)=\frgo\oplus\frguno$,
where $\frgo=\inder(T)$ and $\frguno=V\otimes T$, with
$\bigl(V,\langle.\vert.\rangle\bigr)$ is a two dimensional vector
space endowed with a nonzero alternating bilinear form. The
simplicity of $T$ shows that $\frguno$ is a sum of two copies of
an irreducible module for $\frgo$. A close look at the
classification in \cite{Kac} reveals that the only possibility for
$\frg(T)$ is to be isomorphic to $\frpsl_{2,2}(k)$, so that
$\inder(T)\cong
\frsl_2(k)\oplus\frsl_2(k)\cong\frsp(V)\oplus\frsp(V)$ and $T$ is
the four dimensional irreducible module $V\otimes V$ (the
$i^{\text{th}}$ copy of $\frsp(V)$ acts on the $i^{\text{th}}$
slot of $V\otimes V$, $i=1,2$). Therefore, $\dim_kT=4$ and $T$ is
endowed with a nondegenerate symmetric bilinear form $(.\vert .)$,
invariant under $\inder(T)$, as so does $V\otimes V$ with respect
to $\frsp(V)\oplus\frsp(V)$ (the symmetric bilinear form is just
the tensor product of the alternating forms on each copy of $V$).
Hence the formula $\Phi(x,y,z,t)=([xyz]\vert t)$ defines a
multilinear skew-symmetric map and the proof is complete.
\end{proof}

\medskip

There appears the natural open question of the classification of
the simple (null) \OTSs\ over fields of prime characteristic. All
the examples in Theorems \ref{th:OTSchar0} and
\ref{th:NullOTSchar0} have their counterparts in prime
characteristic, with only a couple of changes needed for the
simple \OTSs\ of G and F types in characteristic $3$. The F-type
\OTS\ becomes a simple null \OTS\ in characteristic $3$. On the
other hand, given a Cayley algebra over a field $k$ of
characteristic $3$, the linear span of the derivations $D_{x,y}$'s
form the unique simple ideal $\inder(C)$, of dimension $7$, of the
fourteen dimensional Lie algebra of derivations $\der(C)$, as
shown in \cite{Albercaetal}, where it is proved that $\inder(C)$
is a form of the simple Lie algebra $\frpsl_3(k)$ (the simple Lie
algebra of type $A_2$ in characteristic $3$). Moreover, any form
of $\frpsl_3(k)$ is the inner derivation algebra of a Cayley
algebra. Hence, over fields of characteristic $3$, the simple
\OTSs\  of G-type make sense, but their inner derivation algebras
are seven dimensional, instead of fourteen dimensional, and thus
$\dim\frg(T)=(3+7)+(2\times 7)=24$.

The proofs of Theorems \ref{th:OTSchar0} and \ref{th:NullOTSchar0}
are based in the classification by Kac of the simple finite
dimensional Lie superalgebras over algebraically closed field of
characteristic $0$, so different methods are needed in prime
characteristic.

Here, we will content ourselves with showing that over fields of
characteristic $3$, there is a new family of simple \OTSs\ related
to Jordan algebras of nondegenerate cubic forms with basepoint.

\begin{examples}\label{ex:Jordanchar3}
Let $J=\calJord(n,1)$ be the Jordan algebra of a nondegenerate
cubic form $n$ with basepoint $1$, over a field $k$ of
characteristic $3$, and assume that $\dim_kJ\geq 3$. Then any
$x\in J$ satisfies a cubic equation \cite[II.4.2]{McCrimmon}
\begin{equation}\label{eq:degree3}
x^{\cdot 3}-t(x)x^{\cdot 2}+s(x)x-n(x)1=0,
\end{equation}
where $t$ is its trace linear form, $s(x)=t(x^\sharp)$ is the spur
quadratic form and the multiplication in $J$ is denoted by $x\cdot
y$.

Let $J_0=\{x\in J:t(x)=0\}$ be the subspace of zero trace
elements. Since $\charac k=3$, $t(1)=0$, so that $k\, 1\in J_0$.
Consider the quotient space $\hat J=J_0/k\, 1$. For any $x\in
J_0$, $s(x)=-\tfrac{1}{2}t(x^{\cdot 2})$ and, by linearization of
\eqref{eq:degree3}, we get that for any $x,y\in J_0$:
\begin{equation}\label{eq:yyx}
\begin{split}
y^{\cdot 2}\cdot x-(x\cdot y)\cdot y&\equiv
 -2t(x,y)y-t(y,y)x \mod k\,1,\\
 &\equiv t(x,y)y-t(y,y)x \mod k\, 1.
\end{split}
\end{equation}
Let us denote by $\hat x$ the class of $x$ modulo $k\,1$. Since
$J_0$ is the orthogonal of $k\,1$ relative to the trace bilinear
form $t(a,b)=t(a\cdot b)$, $t$ induces a nondegenerate symmetric
bilinear form on $\hat J$ defined by $t(\hat x,\hat y)=t(x,y)$ for
any $x,y\in J_0$. Now, consider for any $x,y\in J_0$ the inner
derivation of $J$ given by $D_{x,y}:z\mapsto x\cdot (y\cdot
z)-y\cdot(x\cdot z)$ (see \cite{JacobsonJordan}). Since the trace
form is invariant under the Lie algebra of derivations, $D_{x,y}$
leaves $J_0$ invariant, and obviously satisfies $D_{x,y}(1)=0$, so
it induces a map $d_{x,y}:\hat J\rightarrow \hat J$, $\hat
z\mapsto \widehat{D_{x,y}(z)}$ and a well defined bilinear map
$\hat J\times \hat J\rightarrow \frgl(\hat J)$, $(\hat x,\hat
y)\mapsto d_{x,y}$. Consider now the triple  product $[...]$ on
$\hat J$ defined by
\[
[\hat x\hat y\hat z]=d_{x,y}(\hat z)
\]
for any $x,y,z\in J_0$. This is well defined and satisfies
\eqref{eq:OTSa}, because of the skew-symmetry of $d_{.,.}$. Also,
\eqref{eq:yyx} implies that
\[
\begin{split}
[\hat x\hat y\hat y]=d_{x,y}(\hat y)&=t(x,y)\hat y-t(y,y)\hat x\\
&=t(\hat x,\hat y)\hat y-t(\hat y,\hat y)\hat x,
\end{split}
\]
so that \eqref{eq:OTSb} is satisfied too, relative to the trace
bilinear form. Since $D_{x,y}$ is a derivation of $J$ for any
$x,y\in J$, \eqref{eq:OTSc} follows immediately, while
\eqref{eq:OTSd} is a consequence of the previous ones.

Therefore, by the nondegeneracy of the trace form
\begin{center}
\emph{$\bigl(\hat J,[...], t(.,.)\bigr)$ is a simple \OTS\ over
$k$.}
\end{center}
\end{examples}

\medskip

There are two possibilities for such Jordan algebras, either
$J=k\times\calJord(q,e)$, where $\calJord(q,e)$ is the Jordan
algebra of a nondegenerate quadratic from $q$ with basepoint $e$,
or $J$ is a central simple Jordan algebra of degree $3$.

In the first case, $\calJord(q,e)=k\,e\oplus W$, where $W=\{x\in
\calJord(q,e): q(x,e)=0\}$, and $\hat J=J_0/k\,1$ can be
identified with $W$. Moreover, for any $x,y,z\in W$ (see
\cite[II.3.3]{McCrimmon}):
\[
\begin{split}
x\cdot (y\cdot z)-y\cdot (x\cdot z)&=
 x\cdot\left(-\tfrac{1}{2}q(y,z)e\right)-y\cdot\left(-\tfrac{1}{2}q(x,z)e\right)\\
 &=q(y,z)x-q(x,z)y,
\end{split}
\]
and
\[
t(x,y)=t(x\cdot y)=-\tfrac{1}{2}t\bigl(q(x,y)e\bigr)=-q(x,y).
\]
Therefore, $[xyz]=t(x,z)y-t(y,z)x$ and $\bigl(\hat
J,[...],t(.,.)\bigr)$ is a simple \OTS\ of orthogonal type. Hence
nothing new appears in this case.

\smallskip

In the second case (\emph{Jordan type}), $J$ is a central simple
Jordan algebra of degree $3$ so, after a scalar extension, $J$ is
the algebra of hermitian matrices $H_3(C)$, where $C$ is either
$k$, $k\times k$, $\Mat_2(k)$ or the algebra of split octonions.
Besides, according to \cite[Theorem VI.9 and Theorem
IX.17]{JacobsonJordan}, the Lie algebra of derivations of $J$,
$\der(J)$, is isomorphic, respectively, to $\frso_3(k)$,
$\frpgl_3(k)$, $\frsp_6(k)$, or the $52$ dimensional simple Lie
algebra of type $F_4$. Also, since $J_0$ is an invariant subspace
for $\der(J)$, $D_{J_0,J_0}$ is an ideal of $\der(J)$. Hence, by
simplicity, $\inder(\hat J)$ is isomorphic either to $\frso_3(k)$,
$\frsp_6(k)$ or the Lie algebra of type $F_4$ if $\dim C=1,4$ or
$8$. Finally, if $\dim C=2$, then $J\cong \Mat_3(k)^+$, where
$x\cdot y=\tfrac{1}{2}(xy+yx)$ for any $x,y\in \Mat_3(k)$, so
\[
D_{x,y}:z\mapsto
\tfrac{1}{4}\Bigl(x(yz+zy)+(yz+zy)x-(y(xz+zx)-(xz+zx)y\Bigr)=
\tfrac{1}{4}[[x,y],z],
\]
for any $x,y$ and, hence, $D_{J_0,J_0}$ is isomorphic to the seven
dimensional simple Lie algebra $\frpsl_3(k)$, and so is
$\inder(\hat J)$.

Therefore, for $J=H_3(k)$, $\dim\hat J=4$ and $t([\hat x\hat y\hat
z],\hat u)$ gives a skew-symmetric nonzero four-linear map
(because the trace is $\der(J)$-invariant). Hence $\bigl(\hat
J,[...],t(.,.)\bigr)$ is a simple \OTS\ of D$_1$-type, as the
dimension of $\inder\hat J\cong\frso_3(k)$ is $3$ (see
\ref{th:OTSchar0}). Again, nothing new appears in this case.

For $J=\Mat_3(k)^+$, up to isomorphism $\hat J=\frpsl_3(k)$ and
$d_{\hat x,\hat y}=\tfrac{1}{4}\ad[\hat x,\hat y]$ for any $\hat
x,\hat y\in \hat J$ (brackets in the Lie algebra $\frpsl_3(k)$).
It follows that $\bigl(\hat J,[...],t(.,.)\bigr)$ is a simple
\OTS\ of G-type.

Finally, for $J=H_3(C)$, with $\dim C=4$ or $8$, $\dim\hat J=13$
or $25$ respectively and the situation has no counterpart in
characteristic $0$. \qed

\bigskip

Summarizing some of the previous work, some other simple Lie
superalgebras in characteristic $3$, with no counterpart in
characteristic $0$, have been found:

\begin{theorem}\label{th:othernewchar3}
Let $k$ be an algebraically closed field of characteristic $3$.
Then there are simple finite dimensional Lie superalgebras $\frg$
over $k$ satisfying:
\begin{romanenumerate}
\item
 $\dim\frg=24\,\bigl(=(3+7)+(2\times 7)\bigr)$,
 $\frgo$ is the direct sum of $\frsl_2(k)$ and
$\frpsl_3(k)$ and, as a $\frgo$-module, $\frguno$ is the tensor
product of the natural two dimensional module for $\frsl_2(k)$ and
the adjoint module for $\frpsl_3(k)$ (G-type).
\item
 $\dim\frg=37\,\bigl(=21+(2\times 8)\bigr)$,
 $\frgo=\frso_7(k)$ and, as a $\frgo$-module, $\frguno$ is the
 direct sum of two copies of its spin module (F-type).
\item
 $\dim\frg=50\,\bigl(=(3+21)+(2\times 13)\bigr)$,
 $\frgo$ is the direct sum of $\frsl_2(k)$ and
$\frsp_6(k)$ and, as a $\frgo$-module, $\frguno$ is the tensor
product of the natural two dimensional module for $\frsl_2(k)$ and
of a $13$ dimensional irreducible module for $\frsp_6(k)$ (namely,
the quotient of the space of zero trace symmetric matrices in
$\Mat_6(k)$ relative to the symplectic involution modulo the
scalar matrices).
\item
 $\dim\frg=105\,\bigl(=(3+52)+(2\times 25)\bigr)$,
 $\frgo$ is the direct sum of $\frsl_2(k)$ and
 of the central simple Lie algebra of type $F_4$ and, as a $\frgo$-module,
 $\frguno$ is the tensor
product of the natural two dimensional module for $\frsl_2(k)$ and
a $25$ dimensional irreducible module for $F_4$.
\end{romanenumerate}
\end{theorem}

\section{Lie algebras and \OTSs}\label{se:OTSLie}

Characteristic $3$ is special too for \OTSs, because it turns out
that any such system is a Lie triple system, that is, the odd part
of a $\bZ_2$-graded Lie algebra. More precisely:

\begin{theorem}\label{th:OTSLie}
Let $\T$ be either an \OTS\ or a null \OTS\ over a field of
characteristic $3$. Define the $\bZ_2$-graded algebra
$\frtg=\frtg(T)=\frtg\subo\oplus\frtg\subuno$, with:
\[
\frtg\subo=\inder(T),\qquad\frtg\subuno=T,
\]
and anticommutative multiplication given by:
\begin{itemize}
\item
$\frtg\subo$ is a Lie subalgebra of $\frtg$;
\item
$\frtg\subo$ acts naturally on $\frtg\subuno$, that is,
$[d,x]=d(x)$ for any $d\in\inder(T)$ and $x\in T$;
\item
$[x,y]=d_{x,y}=[xy.]$, for any $x,y\in T$.
\end{itemize}
Then $\frtg(T)$ is a Lie algebra. Moreover, $T$ is simple if and
only if so is $\frtg(T)$.
\end{theorem}
\begin{proof}
As in the proof of Theorem \ref{th:STSsuperLie}, if $\T$ is an
\OTS\ and $x,y,z\in T$, in $\tilde\frg(T)$:
\[
\begin{split}
[[x,y],z]+&[[y,z],x]+[[z,x],y]\\
  &=[xyz]+[yzx]+[zxy]\\
  &=[xyz]+[yzx]-2[zxy]\\
  &=\bigl([xyz]+[xzy]\bigr)-\bigl([zyx]+[zxy]\bigr)\quad\text{(by
  \eqref{eq:OTSa})}\\
  &=\bigl(\sigma_{x,y}(z)+\sigma_{x,z}(y)\bigr) -
    \bigl(\sigma_{z,y}(x)+\sigma_{z,x}(y)\bigr)\quad\text{(by
    \eqref{eq:OTSb})}\\
  &=3\bigl((x\vert y)z-(y\vert z)x\bigr)=0
\end{split}
\]
and the same argument works for null \OTSs\ with trivial
$(.\vert.)$. The rest of the proof follows exactly as in Theorem
\ref{th:STSsuperLie}.
\end{proof}

\medskip

Let $k$ be a field of characteristic $3$ and consider the
analogues over $k$ of the \OTSs\ that appear in Theorems
\ref{th:OTSchar0} and \ref{th:NullOTSchar0}. We finish the paper
by computing the corresponding Lie algebras $\frtg$. Besides some
classical Lie algebras, the simple Lie algebras $L(\epsilon)$ of
Kostrikin, as well as Brown's $29$-dimensional simple Lie algebra,
appear again.

\medskip
\begin{examples}\label{ex:OTSexamples} \null

\noindent\emph{Orthogonal type}:\quad Here $[xyz]=\sigma_{x,y}(z)$
for any $x,y,z\in T$, so that $\inder(T)=\frso(T)$ and
$\frtg(T)=\frso(T)\oplus T$, which is isomorphic to the orthogonal
Lie algebra of a space which is the orthogonal sum of $T$ and a
one dimensional subspace.

\medskip

\noindent\emph{Unitarian type}:\quad Consider the split case:
$K=k\times k$. Hence, as in \cite[p.\,358]{EKO}, we may assume
that $T=W\oplus W^*$ for a vector space $W$ with $\dim_kW\geq 3$
and, after scaling, we get that the triple product is determined
by $[WWT]=0=[W^*W^*T]$ and by
\begin{subequations}\label{eq:OTSxfyg}
\begin{align}
[xfy]&=f(x)y-2f(y)x,\label{eq:OTSxfy}\\
[xfg]&=-f(x)g+2g(x)f,\label{eq:OTSxfg}
\end{align}
\end{subequations}
for any $x,y\in W$ and $f,g\in W^*$ (compare with the
\emph{special case} in Examples \ref{ex:examples}). The
restriction map $\inder(T)\rightarrow \frgl(W)$, $d\mapsto
d\vert_W$ gives an isomorphism between $\inder(T)$ and either
$\frsl(W)$ or $\frgl(W)$, depending on the dimension of $W$ being
or not congruent to $2$ modulo $3$. Let $\tilde W$ be the
$\bZ_2$-graded vector space with $\tilde W\subo=k\, e$ and $\tilde
W\subuno=W$ and identify the special linear Lie algebra
$\frsl(\tilde W)$ with the vector space of $2\times 2$ matrices
$\bigl\{\bigl(\begin{smallmatrix} -\trace A&f\\
x&A\end{smallmatrix}\bigr) :x\in W,\, f\in W^*,\, A\in
\frgl(W)\bigr\}$. Let $\tilde\frz$ be the center of $\frsl(\tilde
W)$ ($\tilde\frz=0$ if $\dim W\not\equiv 2$ modulo $3$, and
$\dim\tilde\frz=1$ otherwise), and let $\frpsl(\tilde
W)=\frsl(\tilde W)/\tilde\frz$ be the corresponding projective
special linear Lie algebra. Then the linear map
$\tilde\Phi:\frtg(T)\rightarrow \frpsl(\tilde W)$, such that:
\[
\begin{split}
\tilde\Phi\bigl((x,f)\bigr)&=\begin{pmatrix} 0&f\\ x&0\end{pmatrix} +\tilde\frz\,,\\
\tilde\Phi(d)&=\begin{cases} \begin{pmatrix} 0&0\\
0&d\vert_W\end{pmatrix}+\tilde\frz\qquad\text{if $\dim W\equiv 2\pmod 3$,}\\
\noalign{\smallskip}
\begin{pmatrix} -\tfrac{\trace d\vert_W}{1+\dim W}&0\\
 0&d\vert_W-\tfrac{\trace d\vert_W}{1+\dim W}id\end{pmatrix}
 +\tilde\frz\qquad\text{otherwise},
\end{cases}
\end{split}
\]
for any $x\in W$, $f\in W^*$ and $d\in \inder(T)$, is an
isomorphism of Lie algebras.

\medskip

\noindent\emph{Symplectic type}:\quad  Consider again the split
case, where $Q=\End_k(V)$, for a two dimensional vector space $V$
endowed with a nonzero alternating bilinear form $\langle .\vert
.\rangle$. Then the arguments in \cite[p.\,359]{EKO} give that, up
to isomorphism, $T=V\otimes W$ for some even dimensional vector
space $W$ of dimension $\geq 4$, endowed with a nondegenerate
alternating bilinear form $\psi:W\times W\rightarrow k$ such that
the triple product becomes (just take $\varphi(u,v)=2\langle
u\vert v\rangle$ in \cite{EKO}):
\[
[(a\otimes x)(b\otimes y)(c\otimes z)]=
\psi(x,y)\gamma_{a,b}(c)\otimes z-2\langle a\vert b\rangle
c\otimes \psi_{x,y}(z)
\]
for any $a,b,c\in V$ and $x,y,z\in W$, where $\gamma_{a,b}=\langle
a\vert .\rangle b+\langle b\vert .\rangle a$ and
$\psi_{x,y}=\psi(x,.)y+\psi(y,.)x$. Besides, $\inder(T)$ is
isomorphic naturally to the direct sum $\frsp(V)\oplus\frsp(W)$ of
the corresponding symplectic Lie algebra. Thus we may identify
$\frtg(T)$ with the $\bZ_2$-graded Lie algebra
\[
\frtg=\bigl(\frsp(V)\oplus\frsp(W)\bigr)\oplus (V\otimes W),
\]
where
\[
[a\otimes x,b\otimes y]=\psi(x,y)\gamma_{a,b}+\langle a\vert
b\rangle \psi_{x,y}\in\frsp(V)\oplus\frsp(W),
\]
for any $a,b\in V$ and $x,y\in W$. Thus $\frtg$ is isomorphic to
the symplectic Lie algebra of the orthogonal sum of $V$ and $W$:
$\frsp(V\perp W)$.

\medskip

\noindent\emph{G-type}:\quad In characteristic $3$, for any $x,y$
in a Cayley algebra $C$:
\[
\begin{split}
D_{x,y}&=[L_x,L_y]+[L_x,R_y]+[R_x,R_y]=[L_x-R_x,L_y-R_y]\\
&=[\ad x,\ad y]=\ad[x,y],
\end{split}
\]
(because $[L_x,R_x]=0$ for any $x$ and
$[[x,y],z]+[[y,z],x]+[[z,x],y]=6\bigl((xy)z-x(yz)\bigr)=0$
\cite[proof of Theorem 1]{Albercaetal}). Therefore, the triple
product on $T=C_0$ in \eqref{eq:OTSG} becomes
\[
[xyz]=\alpha[[x,y],z],
\]
and thus $\frtg(T)=\inder(T)\oplus T\cong
C_0\otimes_kk[\epsilon]$, where $k[\epsilon]=k1\oplus k\epsilon$,
with $\epsilon^2=\alpha$, and $C_0$ is a Lie algebra with the
usual bracket $[x,y]=xy-yx$, which is a form of $\frpsl_3(k)$ (the
split simple Lie algebra of type $A_2$ over $k$). Note that if
$\alpha\in k^2$, $\frtg(T)\cong C_0\oplus C_0$.

\medskip

\noindent\emph{F-type}:\quad Here $\frtg(T)$ is a $\bZ_2$-graded
Lie algebra whose even part is a simple Lie algebra of type $B_3$
and the even part is its spin module. Therefore $\frtg(T)$ is a
form of the $29$ dimensional simple Lie algebra discovered by
Brown \cite{Brown29}.

\medskip

\noindent\emph{D$_\mu$-type}:\quad Assume here that the ground
field $k$ is algebraically closed and consider the simple Lie
superalgebra $\frg=D(2,1;\alpha)=\Gamma(1,\alpha,-(1+\alpha))$,
$\alpha\ne 0,-1$ (notation as in \cite[pp.\,17-18]{Scheunert}):
\[
\begin{cases}
\frgo=\frsp(V_1)\oplus\frsp(V_2)\oplus\frsp(V_3),\\
\frguno=V_1\otimes V_2\otimes V_3,
\end{cases}
\]
where the $V_i$'s are two dimensional vector spaces endowed with
nonzero alternating bilinear forms $\langle.\vert.\rangle$ (the
same notation will be used for the three $V_i$'s), and the Lie
bracket of even or even and odd elements is the natural one, while
\begin{multline*}
[x_1\otimes x_2\otimes x_3,y_1\otimes y_2\otimes y_3]\\
=\langle x_2\vert y_2\rangle\langle x_3\vert
y_3\rangle\gamma_{x_1,y_1} +
 \alpha \langle x_1\vert y_1\rangle\langle x_3\vert y_3\rangle\gamma_{x_2,y_2}
 -(1+\alpha)\langle x_1\vert y_1\rangle\langle x_2\vert
y_2\rangle\gamma_{x_3,y_3},
\end{multline*}
for any $x_i,y_i\in V_i$, $i=1,2,3$. According to Theorem
\ref{th:OTSsuperLie}, $T_\alpha=V_2\otimes V_3$ is a simple \OTS\
with the triple product and symmetric bilinear form determined by
\[
\begin{split}
&[(x_2\otimes x_3)(y_2\otimes y_3)(z_2\otimes z_3)]\\
&\qquad\qquad =\alpha \langle x_3\vert
y_3\rangle\gamma_{x_2,y_2}(z_2)\otimes z_3
 -(1+\alpha)\langle x_2\vert
y_2\rangle z_2\otimes\gamma_{x_3,y_3}(z_3),\\
 \noalign{\smallskip}
&\bigl( x_2\otimes x_3\vert y_2\otimes y_3\bigr)=\langle x_2\vert
y_2\rangle\langle x_3\vert y_3\rangle,
\end{split}
\]
for arbitrary elements $x_i,y_i,z_i\in V_i$, $i=2,3$. These are
precisely the \OTSs\ of D$_\mu$-type (see \cite{EKO}) with $\mu\ne
1$. Here the simple $\bZ_2$-graded Lie algebra $\frtg(T_\alpha)$
satisfies
\[
\begin{cases}
\frtg\subo=\frsp(V_2)\oplus\frsp(V_3),\\
\frtg\subuno= V_2\otimes V_3,
\end{cases}
\]
and the bracket of odd basis elements is given by:
\[
[ x_2\otimes x_3, y_2\otimes y_3] =
 \alpha \langle x_3\vert y_3\rangle\gamma_{x_2,y_2}
 -(1+\alpha)\langle x_2\vert y_2\rangle\gamma_{x_3,y_3},
\]
for $x_i,y_i\in V_i$, $i=2,3$. Proposition \ref{pr:Lepsilon} shows
that $\frtg(T_\alpha)$ is isomorphic to the Kostrikin Lie algebra
$L\left(\tfrac{1+\alpha}{\alpha}\right)$. Thus we get all the
Kostrikin algebras $L(\epsilon)$ for $\epsilon\ne 1$ ($L(1)$ will
be obtained shortly). However, for $\alpha=-1$, we get a Lie
superalgebra $\frg=\frgo\oplus\frguno$ with
\begin{equation}\label{eq:Luno}
\begin{cases}
\frgo=\frsp(V_1)\oplus\frsp(V_2),\\
\frguno=V_1\otimes V_2\otimes V_3,
\end{cases}
\end{equation}
(the $V_i$'s as before) and Lie bracket determined
by
\begin{multline}\label{eq:moreLuno}
[x_1\otimes x_2\otimes x_3,y_1\otimes y_2\otimes y_3]\\
=\langle x_2\vert y_2\rangle\langle x_3\vert
y_3\rangle\gamma_{x_1,y_1}
 - \langle x_1\vert y_1\rangle\langle x_3\vert y_3\rangle\gamma_{x_2,y_2},
\end{multline}
for $x_i,y_i\in V_i$, $i=1,2,3$. This Lie superalgebra is
isomorphic to $\frpsl_{2,2}(k)$ and shows that $T_{-1}=V_2\otimes
V_3$ is a simple \OTS\ with the triple product and symmetric
bilinear form determined by
\[
\begin{split}
&[(x_2\otimes x_3)(y_2\otimes y_3)(z_2\otimes z_3)] =- \langle
x_3\vert y_3\rangle\gamma_{x_2,y_2}(z_2)\otimes
z_3,\\
&\bigl( x_2\otimes x_3\vert y_2\otimes y_3\bigr)=\langle x_2\vert
y_2\rangle\langle x_3\vert y_3\rangle,
\end{split}
\]
for arbitrary elements $x_i,y_i,z_i\in V_i$, $i=2,3$. From
\cite[Proposition 3.3]{EKO} it follows that  $T_{-1}$ corresponds
to the D$_1$-type. The associated $\bZ_2$-graded Lie algebra
$\frtg=\frtg(T_{-1})$ satisfies that $\frtg\subo=\frsp(V_2)$ and
$\frtg\subuno=V_2\otimes V_3$ (the direct sum of two copies of the
natural module for $\frtg\subo$), and this is isomorphic to
$\frpsl_3(k)$. Also, because of \eqref{eq:Luno} and
\eqref{eq:moreLuno}, it turns out that $T=V_1\otimes V_2$, with
triple product $[xyz]=d_{x,y}(z)$ determined by
\[
d_{x_1\otimes x_2,y_1\otimes y_2}=\langle x_2\vert y_2\rangle
\gamma_{x_1,y_1}-\langle x_1\vert y_1\rangle \gamma_{x_2,y_2}
\]
for $x_i,y_i\in V_i$, $i=1,2$, is a simple null \OTS, and its
associated Lie algebra $\frtg(T)$ is, because of Proposition
\ref{pr:Lepsilon}, isomorphic to the Kostrikin algebra $L(1)$.

\medskip

\noindent\emph{Jordan type}:\quad Let $J$ be a central simple
Jordan algebra of degree $3$, so that $\dim_kJ=6,9,15$ or $17$,
and let $\bigl(\hat J,[...],t(.,.)\bigr)$ be its associated \OTS\
of Jordan type as in Examples \ref{ex:Jordanchar3}. Then the
arguments given there show that $\der(J)=[der(J),der(J)]\simeq
\inder(\hat J)$ if $\dim_kJ\ne 9$, and
$[\der(J),\der(J)]=\inder(\hat J)$ if $\dim_kJ=9$. Associated to
$J$ there is the Lie algebra $\frL_0(J)=L_{J_0}\oplus D_{J,J}$ (a
subalgebra of $\frgl(J)$), where $L_x$ denotes the left
multiplication by $x$ (see \cite[VI.9]{JacobsonJordan}).
(Incidentally, this is the Lie algebra that appears in the second
row in Tits construction of Freudenthal's Magic Square if the
characteristic were $\ne 3$ \cite[Theorem 4.13]{Schafer}.) The
derived algebra is $\frL_0'(J)=[\frL_0(J),\frL_0(J)]=L_{J_0}\oplus
D_{J_0,J_0}$, because $D_{J,J}(J_0)=J_0$, and its center is
$k\,L_1$. Then the natural map $\frL_0'(J)\rightarrow \frtg(\hat
J)=\inder(\hat J)\oplus \hat J$, which is the identity on
$D_{J_0,J_0}=\inder(\hat J)$ and takes $L_x$, for $x\in J_0$, to
$\hat x$, induces an isomorphism $\frtg(\hat J)\cong
\frL_0'(J)/k\, L_1$.

If $\dim_kJ=6$, $\hat J$ is of D$_1$-type and $\frtg(\hat J)$ is a
form of $\frpsl_3(k)$, while for $\dim_kJ=9$, $\hat J$ is of
G-type and $\frtg(\hat J)$ is a form of
$\frpsl_3(k)\oplus\frpsl_3(k)$. If $\dim_kJ=15$ and $k$ is
algebraically closed, $J$ is, up to isomorphism, the algebra of
hermitian matrices in $\Mat_6(k)$ relative to the standard
symplectic involution, $D_{J,J}=D_{J_0,J_0}=\inder(\hat J)$ is
given by the adjoint action on $J$ of the skew-hermitian matrices:
$\frsp_6(k)$, and then $\frL_0'(J)\cong \frsl_6(k)$ and hence
$\frtg(\hat J)$ is isomorphic to $\frpsl_6(k)$. Therefore, in
general, if $\dim_kJ=15$, $\frtg(\hat J)$ is a form of
$\frpsl_6(k)$. Finally, if $J$ is exceptional ($\dim_kJ=27$) and
$k$ is algebraically closed, $\frL_0'(J)$ is the Lie algebra of
type $E_6$ (which has a one dimensional center), and $\frtg(\hat
J)$ is then isomorphic to the simple finite dimensional
contragredient Lie algebra $E_6'$ (notation as in \cite[\S
3]{VK}). \qed

\end{examples}


\begin{thebibliography}{BEMN02}

\bibitem[All79]{Allison}
B.~N. Allison, \emph{Models of isotropic simple {L}ie algebras},
Comm. Algebra
  \textbf{7} (1979), no.~17, 1835--1875. 

\bibitem[BDE03]{BDE}
Pilar Benito, Cristina Draper, and Alberto Elduque, \emph{Models
of the
  octonions and {$G\sb 2$}}, Linear Algebra Appl. \textbf{371} (2003),
  333--359. 

\bibitem[BEMN02]{Albercaetal}
Pablo~Alberca Bjerregaard, Alberto Elduque, C{\'a}ndido
  Mart{\'{\i}}n~Gonz{\'a}lez, and Francisco~Jos{\'e} Navarro~M{\'a}rquez,
  \emph{On the {C}artan-{J}acobson theorem}, J. Algebra \textbf{250} (2002),
  no.~2, 397--407. 

\bibitem[Bou02]{Bou}
Nicolas Bourbaki, \emph{Lie groups and {L}ie algebras. {C}hapters
4--6}, Elements of Mathematics (Berlin),
      translated from the 1968 French original by Andrew Pressley,
 Springer-Verlag, Berlin,
      2002.

\bibitem[Bro82]{Brown29}
Gordon Brown, \emph{Properties of a {$29$}-dimensional simple
{L}ie algebra of
  characteristic three}, Math. Ann. \textbf{261} (1982), no.~4, 487--492.

\bibitem[Bro84]{BrownFTS}
\bysame, \emph{Freudenthal triple systems of characteristic
three}, Algebras
  Groups Geom. \textbf{1} (1984), no.~4, 399--441. 

\bibitem[Bro90]{BrownCont}
\bysame, \emph{Structure of certain simple {L}ie algebras of
characteristic
  three}, Lie algebra and related topics (Madison, WI, 1988), Contemp. Math.,
  vol. 110, Amer. Math. Soc., Providence, RI, 1990, pp.~27--31.

\bibitem[Dix84]{Dixmier}
J.~Dixmier, \emph{Certaines alg\`ebres non associatives simples
d\'efinies par
  la transvection des formes binaires}, J. Reine Angew. Math. \textbf{346}
  (1984), 110--128. 

\bibitem[Eld96]{Elduque3fold}
Alberto Elduque, \emph{On a class of ternary composition
algebras}, J. Korean
  Math. Soc. \textbf{33} (1996), no.~1, 183--203. 

\bibitem[Eld04]{MSII}
\bysame, \emph{The Magic Square and Symmetric Compositions II},
preprint.

\bibitem[EO02]{EO}
Alberto Elduque and Susumu Okubo, \emph{Composition
superalgebras}, Comm. Algebra \textbf{30} (2002), no.~11,
      5447--5471.


\bibitem[EKO03]{EKO}
Alberto Elduque, Noriaki Kamiya, and Susumu Okubo, \emph{Simple
{$(-1,-1)$}
  balanced {F}reudenthal {K}antor triple systems}, Glasg. Math. J. \textbf{45}
  (2003), no.~2, 353--372. 


\bibitem[Fau71]{Faulkner}
John~R. Faulkner, \emph{A construction of {L}ie algebras from a
class of
  ternary algebras}, Trans. Amer. Math. Soc. \textbf{155} (1971), 397--408.

\bibitem[Fer72]{Ferrar}
J.~C. Ferrar, \emph{Strictly regular elements in {F}reudenthal
triple systems},
  Trans. Amer. Math. Soc. \textbf{174} (1972), 313--331 (1973).

\bibitem[FF72]{FaulknerFerrar}
John~R. Faulkner and Joseph~C. Ferrar, \emph{On the structure of
symplectic
  ternary algebras}, Nederl. Akad. Wetensch. Proc. Ser. A {\bf 75}=Indag. Math.
  \textbf{34} (1972), 247--256. 


\bibitem[Fre54]{FrII}
\bysame, \emph{Beziehungen der {${\mathfrak E}\sb 7$} und
{${\mathfrak E}\sb 8$} zur
  {O}ktavenebene. {II}}, Nederl. Akad. Wetensch. Proc. Ser. A. \textbf{57}
   = Indag. Math. \textbf{16} (1954), 363--368.

\bibitem[Fre59]{FrVIII}
\bysame, \emph{Beziehungen der {${\mathfrak E}\sb 7$} und
{${\mathfrak E}\sb 8$} zur
  {O}ktavenebene. {VIII}}, Nederl. Akad. Wetensch. Proc. Ser. A. \textbf{62}
   = Indag. Math. \textbf{21} (1959), 447--465.


\bibitem[Jac66]{JacobsonJordan}
N.~Jacobson, \emph{Structure theory for a class of {J}ordan
algebras}, Proc.
  Nat. Acad. Sci. U.S.A. \textbf{55} (1966), 243--251.

\bibitem[Kac77]{Kac}
V.~G. Kac, \emph{Lie superalgebras}, Advances in Math. \textbf{26}
(1977),
  no.~1, 8--96. 

\bibitem[Kan73]{Kantor}
I.~L. Kantor, \emph{Models of the exceptional {L}ie algebras},
Dokl. Akad. Nauk
  SSSR \textbf{208} (1973), 1276--1279. \MR{MR0349779 (50 \#2272)}

  \bibitem[Kan90]{KantorFr}
\bysame, \emph{On the definition of a {F}reudenthal trilinear
operation},
  Sibirsk. Mat. Zh. \textbf{31} (1990), no.~4, 60--67, 221.

\bibitem[KO03]{KamiyaOkubo}
Noriaki Kamiya and Susumu Okubo, \emph{Construction of {L}ie
superalgebras
  {$D(2,1;\alpha),\ G(3)$} and {$F(4)$} from some triple systems}, Proc. Edinb.
  Math. Soc. (2) \textbf{46} (2003), no.~1, 87--98.

\bibitem[Kos70]{Kostrikin}
A.~I. Kostrikin, \emph{A parametric family of simple {L}ie
algebras}, Izv.
  Akad. Nauk SSSR Ser. Mat. \textbf{34} (1970), 744--756.

\bibitem[Kru04]{Krutelevich}
S.~Krutelevich, \emph{Jordan algebras, exceptional groups, and
higher composition laws}, arXiv:math.NT/0411104.

\bibitem[Loo75]{Loos}
Ottmar Loos, \emph{Jordan pairs}, Springer-Verlag, Berlin, 1975.

\bibitem[McC04]{McCrimmon}
Kevin McCrimmon, \emph{A taste of {J}ordan algebras},
Universitext,
  Springer-Verlag, New York, 2004. 

\bibitem[Mey68]{Meyberg}
Kurt Meyberg, \emph{Eine {T}heorie der {F}reudenthalschen
{T}ripelsysteme. {I},
  {II}}, Nederl. Akad. Wetensch. Proc. Ser. A 71=Indag. Math. \textbf{30}
  (1968), 162--174, 175--190. 

\bibitem[Oku93]{OkuboI}
Susumu Okubo, \emph{Triple products and {Y}ang-{B}axter equation.
{I}.
  {O}ctonionic and quaternionic triple systems}, J. Math. Phys. \textbf{34}
  (1993), no.~7, 3273--3291. 

\bibitem[Sch95]{Schafer}
Richard~D. Schafer, \emph{An introduction to nonassociative
algebras}, Dover
  Publications Inc., New York, 1995. 

\bibitem[Sch79]{Scheunert}
Manfred Scheunert, \emph{The theory of {L}ie superalgebras},
Lecture Notes in
  Mathematics, vol. 716, Springer, Berlin, 1979. 

\bibitem[VK71]{VK}
B.~Ju. Ve{\u\i}sfe{\u\i}ler and V.~G. Kac, \emph{Exponentials in
{L}ie algebras
  of characteristic {$p$}}, Izv. Akad. Nauk SSSR Ser. Mat. \textbf{35} (1971),
  762--788. 

\bibitem[YA75]{YamAs}
Kiyosi Yamaguti and Hiroshi Asano, \emph{On the {F}reudenthal's
construction of
  exceptional {L}ie algebras}, Proc. Japan Acad. \textbf{51} (1975), no.~4,
  253--258. 

\end{thebibliography}


\providecommand{\bysame}{\leavevmode\hbox
to3em{\hrulefill}\thinspace}
\providecommand{\MR}{\relax\ifhmode\unskip\space\fi MR }
\providecommand{\MRhref}[2]{%
  \href{http://www.ams.org/mathscinet-getitem?mr=#1}{#2}
} \providecommand{\href}[2]{#2}

\end{document}